\documentclass[a4paper, 10pt, final]{article}

\usepackage[all]{xy}
\usepackage{style}
 \usepackage{amsthm, eufrak}
\usepackage{dsfont}

 \usepackage{geometry}

\newcommand{\CH}{\operatorname{CH}\nolimits}

 \title{Niveau and coniveau filtrations \linebreak on cohomology
   groups and Chow groups}

\author{Charles Vial}
\date{}
\begin{document}

\maketitle

\begin{abstract} The Bloch--Beilinson--Murre conjectures predict the
  existence of a descending filtration on Chow groups of smooth
  projective varieties which is functorial with respect to the action
  of correspondences and whose graded parts depend solely on the
  topology -- i.e. the cohomology -- of smooth projective varieties.
  In this paper, given a smooth projective complex variety $X$, we
  wish to explore, at the cost of having to assume general conjectures
  about algebraic cycles, how the coniveau filtration on the
  cohomology of $X$ has an incidence on the Chow groups of $X$.
  However, by keeping such assumptions minimal, we are able to prove
  some of these conjectures either in low-dimensional cases or when a
  variety is known to have small Chow groups. For instance, we give a
  new example of a fourfold of general type with trivial Chow group of
  zero-cycles and we prove Murre's conjectures for threefolds
  dominated by a product of curves, for threefolds rationally
  dominated by the product of three curves, for rationally connected
  fourfolds and for complete intersections of low degree.  The BBM
  conjectures are closely related to Kimura--O'Sullivan's notion of
  finite-dimensionality. Assuming the standard conjectures on
  algebraic cycles the former is known to imply the latter. We show
  that the missing ingredient for finite-dimensionality to imply the
  BBM conjectures is the coincidence of a certain niveau filtration
  with the coniveau filtration on Chow groups.

\noindent {\it 2010 Mathematics Subject Classification 14C15, 14C25, 14C30}
\end{abstract}

\section*{Introduction}

Chow groups and cohomology groups constitute important invariants of
algebraic varieties. While cohomology groups can usually be computed
explicitly, Chow groups are still largely mysterious and computing
them is a difficult problem in algebraic geometry. On the arithmetic
side it is expected via Beilinson's conjecture that the Chow groups of
smooth projective varieties defined over number fields are finite in
some sense. Here we wish to focus on the geometric side. Given a
smooth projective variety $X$ defined over the complex numbers,
Mumford \cite{Mumford} (in the case of surfaces), Roitman
\cite{Roitman} and Bloch and Srinivas \cite{BS} proved that if the
Chow group $\CH_0(X)$ of $0$-cycles on $X$ is generated by a point then
the cohomology groups of $X$ are supported on a divisor in positive
degree. In other words if $\CH_0(X)=\Q$ then $H^i(X) = N^1H^i(X)$ for
all $i>0$. Here $N$ is the so-called coniveau filtration. Bloch has
asked if the converse holds true, that is if the coniveau filtration
on the cohomology of $X$ controls its Chow groups.  More generally as
part of the Bloch--Beilinson philosophy it is expected that the
coniveau filtration on cohomology should reflect on Chow groups.

Rather than working with the coniveau filtration $N$ on cohomology, we
introduce a finer filtration -- the niveau filtration $\N$ on homology
-- whose definition and relation with the usual coniveau filtration
$N$ is given in \S1.1. The relevance of $\N$ can be read off
Proposition \ref{relevance} which shows that if $\CH_0(X)=\Q$ then
$H^i(X) = \N^1H^i(X)$ for all $i>0$. Proposition
\ref{filtration-coincide} shows, however, that these two filtrations
are expected to be the same. This niveau filtration was implicitly
used by Schoen in \cite{Schoen} in order to prove his Theorem 0.3
which we recall below as Theorem \ref{LS}.

In this paper, we want to make more precise the link between Chow
groups and the support of cohomology groups.  A first step consists in
giving an algebraic origin to the niveau filtration on homology. This
is embodied in Theorem \ref{Ksupport}. The main general result of this
first section is then Theorem \ref{action}. In order to state Theorem
\ref{Ksupport} we have to restrict our attention to a certain class of
varieties, those that satisfy the property $(\star)$ below. Likewise
Theorem \ref{action} only applies to those varieties that in addition
to $(\star)$ satisfy the property $(\star \star)$ below.  Let's
immediately mention that those conditions are not empty; abelian
varieties and Fermat varieties \cite{KS} of dimension less than or equal
to $5$ satisfy those conditions. More importantly, Theorems
\ref{Ksupport} and \ref{action} make it possible to exhibit some new
examples of varieties, the Chow groups of which can be computed, see
\S2.3.  Moreover, the properties $(\star)$ and $(\star \star)$ are
expected to be satisfied by all smooth projective varieties by general
conjectures on algebraic cycles.

In the second section, we consider a variety $X$ that satisfies the
properties $(\star)$ and $(\star \star)$ and we derive several
propositions relating the support of the cohomology groups of $X$ to
the Chow groups of $X$. Proposition \ref{gen} gives the most general
statement. Proposition \ref{support} treats the case of zero-cycles
and proves that if there is an integer $i$ such that $H_k(X)=
\widetilde{N}^1H_k(X)$ for all $k > i$, then $\CH_0(X)$ is supported in
dimension $i$. This proposition generalizes previous results : Kimura
treated the case when $X$ is a surface \cite[Cor.  7.7]{Kimura} and
Voisin treated the case $i=0$ assuming Grothendieck's generalized
Hodge conjecture \cite[Th.  3]{Voisin2}. This section culminates in
\S2.3 where we exhibit a new example of a fourfold of general type $Y$
whose Chow group of zero-cycles is generated by a point. Precisely,
$Y$ is the quotient of the Fermat hypersurface $\{x_0^7 + x_1^7 +
\ldots +x_{5}^7 =0\} \subset \P_\C^{5}$ by an appropriate free action
of the group of seventh roots of unity. Given Proposition \ref{gen},
the key point consists in proving Theorem \ref{GHC} which implies that
$H_4(Y) = \N^1H_4(Y)$, in other words, that the cohomology of $Y$ is
generated by that of a surface.

In the spirit of Jannsen's paper \cite{Jannsen4}, we are then
interested in relating various conjectures on algebraic cycles.
Grothendieck's standard conjectures \cite{Kleiman} imply that
homological equivalence of algebraic cycles coincides with the
coarsest of equivalence relations, namely numerical equivalence.
Kimura's finite-dimensionality conjecture \cite{Kimura} (independently
stated by O'Sullivan) is an attempt to fill in the gap there exists
between the finest equivalence relation, i.e. rational equivalence,
and homological equivalence of cycles.  However this conjecture is
less precise, i.e.  weaker (see \cite{AnKa}), than the
Bloch--Beilinson--Murre conjectures \cite{Jannsen} which describe such
a gap in terms of the cohomology of $X$. In \S3, we wish to give a
sufficient condition for Kimura's conjecture to imply Murre's
conjectures (these are equivalent to Bloch's and Beilinson's by
Jannsen \cite{Jannsen}). For this purpose, assuming Grothendieck's
standard conjectures and Kimura's conjecture, we define thanks to
Theorem \ref{action} a filtration on the Chow groups of $X$ that
``lifts'' the coniveau filtration on the cohomology of $X$. We prove
in Proposition \ref{end} that if this filtration on Chow groups
coincides with Jannsen's coniveau filtration on Chow groups
\cite[5.10(b)]{Jannsen4} then the BBM conjectures hold.

Finally, in the fourth section, we settle Murre's conjectures (except
perhaps for the ``independency'' conjecture) in some new cases.
Examples of varieties for which Murre's conjectures are known to hold
true include uniruled $3$-folds \cite{dAMS}, $3$-folds with a nef
tangent bundle \cite{Iyer} and elliptic modular $3$-folds \cite{GM}.
Our new examples include $3$-folds rationally dominated by the product
of three curves, rationally connected $3$-folds, Calabi-Yau $3$-folds,
abelian $3$-folds, rationally connected $4$-folds, $4$-folds admitting
a curve as their base for their maximal rationally connected
fibration, rationally connected $5$-folds with vanishing Hodge number
$h^{3,1}$, and some complete intersections of low degree, e.g. cubic
$6$-folds, quartic $7$-folds and the smooth intersection of two
quadrics in $\P^{10}$.  A known case due to Murre \cite{Murre1} that
we recover through our method is given by the product of a curve with
a surface.

\paragraph{Notations.} Throughout, $k$ denotes a subfield of the field
of complex numbers $\C$ and $X$ denotes a smooth projective variety
over $k$ of pure dimension $d$ given with an ample line bundle $L$. We
set $H_i(X) := H_{i}(X(\C),\Q)$; this group is isomorphic to
$H^{2d-i}(X(\C),\Q)$. The Chow group $\CH_l(X)$ is the group of
algebraic cycles of dimension $l$ on $X$ with rational coefficients
modulo rational equivalence. Given $\sim$ an adequate equivalence
relation \cite[1.3]{Jannsen4}, e.g. algebraic or homological
equivalence, $\CH_l(X)_{\sim}$ denotes those $l$-cycles that are $\sim
0$.

\paragraph{Acknowledgements.} I am grateful to the referee for his
careful reading and for his insightful comments, especially for
suggesting Propositions \ref{Ksupportsection-prop} and
\ref{independent-L} and Lemma \ref{special-iso}. Thanks to
Fran\c{c}ois Charles, Burt Totaro and Claire Voisin for their
interest, comments and encouragement. This work started while I was a
student at the \'Ecole Normale Sup\'erieure and a PhD student at
Trinity College, Cambridge.  It is now supported by a Nevile Research
Fellowship at Magdalene College, Cambridge and an EPSRC Postdoctoral
Fellowship under grant EP/H028870/1. I would like to thank all four
institutions for their support.

\section{Niveau and Coniveau filtrations on cohomology}

\subsection*{Statement of the results}\label{statement}

The homology groups of $X$ are endowed with a coniveau filtration
that was studied, for example, by Bloch and Ogus \cite{BO}, but also by
Jannsen \cite{Jannsen}. There is yet another filtration on $H_i(X)$
which is finer than the coniveau filtration and which is called the
niveau filtration. Such a filtration appears implicitly, for instance,
in \cite{Schoen}. The niveau and the coniveau filtrations are both of
interest and Proposition \ref{filtration-coincide} shows that these
filtrations are expected to coincide. \medskip

The \emph{coniveau filtration} on $H_i(X)$ can be defined as follows :
$${N}^j H_i(X) := \sum \im \big(\Gamma_* : H^{i-2j}(Y) \r H_i(X)
\big),$$ where the sum runs through all smooth projective varieties
$Y$ and through all correspondences $\Gamma \in \CH_{i-j}( Y \times
X)$.  The subgroup ${N}^j H_i(X)$ should be thought of as those
classes in $H_i(X)$ that are supported on a closed subscheme of
dimension $i-j$. \medskip

The \emph{niveau filtration} on $H_i(X)$ can be defined as follows :
 $$\widetilde{N}^j H_i(X) := \sum \im \big(\Gamma_* : H_{i-2j}(Y) \r
 H_i(X) \big),$$ where the sum runs through all smooth projective
 varieties $Y$ and through all correspondences $\Gamma \in \CH^{d-j}( Y
 \times X)$. It can be shown that this sum can be restricted to those
 varieties $Y$ of dimension $i-2j$ in which case the correspondence
 $\Gamma$ can be thought of as a family of $j$-cycles parametrized by
 $Y$. In particular we immediately see that $\N^jH_i(X) \subseteq
 N^jH_i(X)$. Because $H_i(X)$ is a finite-dimensional $\Q$-vector
 space, we can find a smooth projective variety $Z_{i,j}$ of pure
 dimension $i-2j$ and a correspondence $\Gamma_{i,j} \in \CH_{i-j}(
 Z_{i,j} \times X)$ such that $\widetilde{N}^j H_i(X) =
 (\Gamma_{i,j})_* H_{i-2j}(Z_{i,j})$. \medskip

 The ample line bundle $L$ on $X$ defines an embedding $X \subseteq
 \P_k^N$ and hence, for any integer $i \leq d$, a map $L^{d-i} :
 H_{2d-i}(X) \r H_{i}(X)$ given by intersecting $d-i$ times with a
 hyperplane. This map is obviously induced by a correspondence on $X
 \times X$ and is an isomorphism of Hodge structures (Hard Lefschetz
 theorem). Given $i \leq d$, we say that $X$ satisfies the property
 $\mathrm{B}_i$ if the isomorphism $L^{i-d} := (L^{d-i})^{-1} :
 H_{i}(X) \r H_{2d-i}(X)$ is induced by a correspondence. If $X$
 satisfies the property $\mathrm{B}_i$ for every $i$, then we say that
 $X$ satisfies the property B. The filtration $\N$ on $H_i(X)$ splits
 canonically into orthogonal pieces as $H_i(X)=\bigoplus_j
 \Gr^j_{\N}H_i(X)$ with respect to the choice of a polarization. If
 $X$ satisfies property B then the splitting is also induced by the
 bilinear form $Q_i := \langle L^{i-d} -, - \rangle$, see Remark
 \ref{niveausplitting}. Condition B is satisfied by abelian varieties
 \cite{Lieberman, Lieberman2}, curves, surfaces, complete
 intersections and any products, hypersurface sections or finite
 quotients thereof. It is a conjecture of Grothendieck (the Lefschetz
 standard conjecture) that all smooth projective varieties should
 satisfy property B. \medskip

 Kahn, Murre and Pedrini \cite[\S 7]{KMP} gave a decomposition of the
 numerical motive of varieties satisfying property B with respect to a
 certain coniveau filtration. Here we wish to refine their result and
 construct a decomposition modulo homological equivalence. For this
 purpose we will have to consider varieties satisfying the following
 property :
 \begin{center} $X$ satisfies B and, for all $i \leq d$ and all $j
   \geq 1$, \linebreak $(\star)$ \hspace{1.4cm} either there exists
   $Z_{i,j}$ as before satisfying $\mathrm{B}_l$ for all $l \leq
   i-2j-2$ \hspace{1.4cm} \linebreak or $N^{j+1}H_i(X) =
   \N^{j+1}H_i(X)$.
\end{center}
Because the property $\mathrm{B}_1$ holds for all smooth projective
varieties, we see that the property $(\star)$ holds for all smooth
projective varieties of dimension at most $5$ that satisfy property B.
In particular, it holds for curves, surfaces, abelian varieties of
dimension $\leq 5$, complete intersections of dimension $\leq 5$,
uniruled $3$-folds, rationally connected varieties of dimension $\leq
4$, rationally connected $5$-folds with $H^3(X,\Omega_X)=0$
\cite{Vial3}. Property $(\star)$ obviously holds for
all varieties if property B holds for all varieties.

\begin{theorem2} \label{Ksupport} If $X$ satisfies $(\star)$ then, for
  all integers $i$ and $j$, there exists a cycle $\pi_{i,j} \in \CH_d(X
  \times X)$ inducing the orthogonal projection $H_*(X) \r
  \Gr_{\N}^jH_i(X) \r H_*(X)$.  Moreover, $\pi_{i,j}$ can be chosen to
  factor through $Z_{i,j}$, i.e.  $\pi_{i,j} = f \circ g$ for some $g
  \in \CH_{d-j}(X \times Z_{i,j})$ and some $f \in \CH_{i-j}(Z_{i,j}
  \times X)$.
\end{theorem2}

With the notations of Kahn, Murre and Pedrini \cite{KMP}, Theorem
\ref{Ksupport} says that the endomorphism $\pi_{i,j}$ of the motive
$h(X)$ factors through the motive $h(Z_{i,j})(j)$.  Let's mention that
the cycles $\pi_{i,j}$ coincide modulo numerical equivalence with the
ones constructed in \cite[7.7.3]{KMP} if $X$ satisfies $(\star)$. This
is proved in Proposition \ref{Ksupportsection-prop}. It follows that
modulo homological equivalence the correspondences $\pi_{i,j}$ are
central idempotents independent of the choice of $L$, and we deduce in
Proposition \ref{independent-L} that the splitting of the niveau
filtration on $H_i(X)$ is independent of the choice of a polarization,
\medskip

Our aim is to show how such a filtration on cohomology reflects on
Chow groups.  For this purpose it is natural to lift the decomposition
of the diagonal $\Delta_X \in \CH_d(X \times X)$ modulo homological
equivalence given above to a decomposition into a sum of mutually
orthogonal idempotents in $\CH_d(X \times X)$. This is achieved (cf.
\cite[Lemma 3.1]{Jannsen}) by considering varieties $X$ satisfying
\begin{center}
  \hspace{-2cm} $(\star \star)$ \hspace{2cm} $\ker \big( cl : \CH_d(X
  \times X) \r H_{2d}(X\times X)\big)$ is a nilpotent ideal.
\end{center}
Here $\CH_d(X \times X)$ is endowed with a ring structure given by the
composition law. For instance, varieties dominated by a product of
curves satisfy $(\star \star)$. This includes abelian varieties and
Fermat varieties \cite{KS}. More generally, any variety that is finite
dimensional in the sense of Kimura satisfies $(\star \star)$; see
\cite[7.2]{Kimura}.

\begin{theorem2} \label{action} Let $X$ be a smooth projective variety
  that satisfies $(\star)$ and $(\star \star)$.  Then there exist
  mutually orthogonal idempotents $\Pi_{i,j} \in \CH_d(X \times X)$
  which are homologically equivalent to the cycles $\pi_{i,j}$ of
  Theorem \ref{Ksupport} and such that $\Delta_X = \sum \Pi_{i,j} \in
  \CH_d(X\times X)$. For any such choice of idempotents, we have

  \begin{enumerate}
  \item $\Pi_{i,j}$ acts as $0$ on $\CH_l(X)$ if either $l <
    j$ or $l>i-j$.
  \item $\Pi_{i,j}$ acts as $0$ on $\CH_l(X)$ if $l = i-j$ and $i<2l$.
  \item $\Pi_{i,j}$ acts as $0$ on $\CH_l(X)$ if $l+1=i-j$ and $i \leq
    2l$.
  \item $\Pi_{i,j} = 0$ if and only if $\Gr_{\widetilde{N}}^j
    H_{i}(X)=0$.
  \item $\CH_i(X)_{\hom} = \ker \big(\Pi_{2i,i} : \CH_i(X) \r
    \CH_i(X) \big).$
  \item in the case the base field $k$ is algebraically closed, if
    $AJ_i : \CH_i(X)_{\alg} \r J_i(X)\otimes \Q$ is the Abel--Jacobi map
    to Griffith's $i^\mathrm{th}$ intermediate Jacobian tensored with
    $\Q$ restricted to algebraically trivial cycles, then
    $$\ker \big( AJ_i \big) = \ker \big(\Pi_{2i+1,i} :
    \CH_i(X)_{\alg} \r \CH_i(X)_{\alg} \big).$$
\end{enumerate}
\end{theorem2}

Examples of varieties satisfying $(\star)$ and $(\star \star)$ include
curves, abelian varieties of dimension $\leq 5$, Fermat varieties of
dimension $\leq 5$, as well as their blow-ups along smooth curves.
More generally, varieties of dimension $\leq 5$ which are dominated by
a product of curves satisfy $(\star)$ and $(\star \star)$.  In
\cite{Vial3} we show that varieties $X$ for which $AJ_i$ is injective
for all $i < d/2$ satisfy $(\star)$ and $(\star \star)$, this includes
Godeaux surfaces \cite{Godeaux}, Fano threefolds \cite{Kollar} and
hypersurfaces of very small degree \cite{Otwinowska}. Moreover, any
variety dominated by one of the varieties above also satisfies
$(\star)$ and $(\star \star)$. \medskip

The first three points in the theorem are consequences of theorem
\ref{Ksupport} together with standard results for varieties satisfying
$(\star \star)$ (Lemmas \ref{liftproj}, \ref{proj}, \ref{proj2} and
\ref{supportlemma}).  The proofs of the two last points are based on
the construction of explicit mutually orthogonal projectors
${\Pi}_{2i,i}$ and ${\Pi}_{2i+1,i}$ as done in \cite{Vial3}, see
Propositions \ref{kerind} and \ref{kerind2}. However, the construction
in \cite{Vial3} is valid under less restrictive hypotheses than the
ones we are working with. For the sake of completeness, we provide a
direct construction of such idempotents under the assumption that $X$
satisfies B, see Theorem \ref{Chowproj}.

\subsection{Filtrations on cohomology}

Both the niveau and the coniveau filtration define decreasing
filtrations on the cohomology ring of a smooth projective variety $X$
over $k$, each of which satisfies $$H_i(X) = N^0H_i(X) \supseteq
N^1H_i(X) \supseteq \ldots \supseteq N^{\lfloor i/2\rfloor}H_i(X)
\supseteq N^{\lfloor i/2\rfloor+1}H_i(X) = 0.$$ The aim of this
paragraph is to compare such filtrations as well as gather well known
and easy properties thereof.

\paragraph{The coniveau filtration $N$.}  This filtration is more
naturally defined on cohomology. The $j^\mathrm{th}$ filtered piece of
$H^i(X)$ defines a rational
sub-Hodge structure of $H^i(X)$ :

\begin{eqnarray}
  \nonumber N^j H^i(X)  & =  & \sum_{Z \subset X} \ker \big(H^i(X)
  \r H^{i}(X-Z) \big) \\
  \nonumber & =  & \sum_{Z \subset X} \im \big(H^i_Z(X) \r
  H^{i}(X) \big) \\
  \nonumber & = & \sum_{f : Y \r X} \im \big(f_* : H^{i-2j}(Y) \r
  H^i(X) \big),
 \end{eqnarray}
 where $Z$ runs through all subschemes of $X$ of codimension $\geq j$,
 and $f : Y \r X$ runs through all morphisms from a smooth projective
 variety $Y$ of pure dimension $\leq \dim X -j$ to $X$. The last
 equality holds in characteristic zero under resolution of
 singularities by Hironaka.  It holds in positive characteristic (for
 $\ell$-adic cohomology) by de Jong's theorem on alterations; see
 Jannsen \cite[7.7]{Jannsen5} and \cite[6.5]{Jannsen4} for more
 details.  We then transcribe this definition in terms of homology and
 set $$ N^j H_i(X) := \sum_{\Gamma \in \CH_{i-j'}(Y \times X)} \im
 \big(\Gamma_* : H^{i-2j'}(Y) \r H_i(X) \big),$$ where the sum runs
 through all integers $j' \geq j$, all smooth projective varieties $Y$
 and all correspondences $\Gamma \in \CH_{i-j'}(Y \times X)$. By
 considering $\P^{j'-j} \times Y$ instead of $Y$ and $\P^{j'-j} \times
 \Gamma$ instead of $\Gamma$ we see that it is superfluous to consider
 those integers $j' >j$. Therefore $$ N^j H_i(X) = \sum_{\Gamma \in
   \CH_{i-j}(Y \times X)} \im \big(\Gamma_* : H^{i-2j}(Y) \r H_i(X)
 \big),$$ where the sum runs through all smooth projective varieties
 $Y$ and all correspondences $\Gamma \in \CH_{i-j}(Y \times X)$. As
 before, using resolution of singularities, we see that actually
 $$ N^j H_i(X) = \sum_{f : Y \r X} \im \big(f_* : H^{i-2j}(Y) \r
 H_i(X) \big),$$ where the sum runs through all morphisms $f : Y \r X$
 from a smooth projective variety $Y$ of pure dimension $\leq i-j$ to
 $X$.

 \paragraph{The niveau filtration $\widetilde{N}$.} This filtration is
 more naturally defined in homology.  By the weak Lefschetz theorem,
 instead of considering varieties $Y$ of dimension $i-j$, we can
 consider smooth linear sections $Z$ of $Y$ of dimension $i-2j$. Such
 a section induces a surjection $H_{i-2j}(Z) \twoheadrightarrow
 H_{i-2j}(Y) \stackrel{\simeq}{\longrightarrow} H_i(Y)$ which is
 induced by an algebraic correspondence if property B holds for $Y$.

 We are thus led to consider families $\Gamma \in \CH_{i-j}(Z \times
 X)$ of $j$-cycles parametrized by $Z$ and look at the image of
 $\Gamma_* : H_{i-2j}(Z) \rightarrow H_{i}(X)$. We then define
 $$\widetilde{N}^j H_i(X) := \sum_{\Gamma \in \CH^{d-j'}( Z \times
   X)} \im \big(\Gamma_* : H_{i-2j'}(Z) \r H_i(X) \big),$$ where the
 union runs through all integers $j'\geq j$, all smooth projective
 varieties $Z$ and all correspondences $\Gamma \in \CH^{d-j'}( Z \times
 X)$. Here we let $j'$ and the dimension of $Z$ vary for greater
 flexibility. However, this is not needed. Indeed, as in the case of
 the coniveau filtration, if $j'>j$, then replace $Z$ by
 $\P^{j'-j}\times Z$ and $\Gamma$ by $\{0\} \times \Gamma \in
 \CH^{d-j}(\P^{j'-j} \times Z \times X)$. Therefore, we can fix $j'$ to
 be equal to $j$ in the sum defining the niveau filtration. Let's now
 show that it is possible to restrict the sum to varieties $Z$ of
 dimension $i-2j$. If $\dim Z > i-2j$, then any smooth linear section
 $\iota :Z' \hookrightarrow Z$ of dimension $i-2j$ induces a
 surjection $H_{i-2j}(Z') \r H_{i-2j}(Z)$.  Replace $Z$ by $Z'$ and
 $\Gamma$ by $\Gamma \circ \iota$. If $\dim Z < i-2j$, then replace
 $Z$ with $\P^a \times Z$ $(a = i-2j - \dim Z)$ and $\Gamma$ by
 $\P^a\times \Gamma \in \CH^{d-j}(\P^a \times Z \times X)$.

 Therefore, we have
 $$\widetilde{N}^j H_i(X) = \sum_{\Gamma\in \CH_{i-j}( Z \times
   X)} \im \big(\Gamma_* : H_{i-2j}(Z) \r H_i(X) \big),$$ where the
 union runs through all smooth projective varieties $Z$ of dimension
 $i-2j$ and all correspondences $\Gamma \in \CH_{i-j}( Z \times X) \
 (=\CH^{d-j}( Z \times X))$.

\paragraph{Properties.}

The last given characterization of the niveau filtration shows
$$\widetilde{N} \subseteq N.$$ Because the Lefschetz isomorphism
$H_i(Z) \r H^i(Z)$ is induced by a correspondence for $i=0$ or $i=1$,
we always have $ \widetilde{N}^{\lfloor i/2\rfloor} H_i(X) =
N^{\lfloor i/2\rfloor} H_i(X).$ More generally, it is expected that
these two filtrations agree:

\begin{proposition} \label{filtration-coincide} Suppose property B
  holds for all smooth projective varieties of dimension $< \dim X$.
  Then, property B holds for $X$ if and only if $N^jH_i(X) =
  \widetilde{N}^jH_i(X)$ for all $i$ and all $j$.
\end{proposition}
\begin{proof}
  Assume $f : Y \r X$ is a morphism of smooth projective varieties
  with $\dim Y = i-j$ and let $Z$ be a smooth linear section of $Y$ of
  dimension $i-2j$. The discussion above shows that the composition
  $H_{i-2j}(Z) \twoheadrightarrow H_{i-2j}(Y)
  \stackrel{\simeq}{\longrightarrow} H_i(Y)
  \stackrel{f_*}{\longrightarrow} H_i(X)$ is induced by a
  correspondence $\Gamma\in \CH_{i-j}( Z \times X)$; it clearly maps
  $H_{i-2j}(Z)$ onto the image of $f_*$, which proves that $N^jH_i(X) =
  \widetilde{N}^jH_i(X)$.

  If $N$ and $\N$ agree on $H_i(X)$ for $i>d$, then in particular
  $N^{i-d}H_i(X)=\N^{i-d}H_i(X)$. It is easy to see from the
  definition that $N^{i-d}H_i(X)= H_i(X)$. Therefore, there is a
  smooth projective variety $Z$ of dimension $2d-i$ and a
  correspondence $\Gamma \in \CH_{d}(Z \times X) $ such that $\Gamma_*
  : H_{2d-i}(Z) \r H_i(X)$ is surjective.  Because $Z$ is assumed to
  satisfy property B, there is a correspondence $s \in \CH_{2d-i}(Z
  \times Z)$ such that $ H_{2d-i}(Z)$ endowed with the pairing
  $\langle - , s_*- \rangle$ is polarized, see \cite{Kleiman}. Hence,
  thanks to Lemma \ref{image} below, the correspondence $\Gamma \circ
  s \circ {}^t\Gamma$ induces an isomorphism $H_{2d-i}(X)
  \stackrel{\simeq}{\longrightarrow} H_i(X)$.  Thus $\alpha := \Gamma
  \circ s \circ {}^t\Gamma \circ L^{i-d}$ induces an automorphism of
  $H_i(X)$. By the theorem of Cayley--Hamilton, its inverse is given
  by $P(\alpha)_*$ for some rational polynomial $P$. It follows that $
  P(\alpha) \circ \Gamma \circ s \circ {}^t\Gamma$ induces the inverse
  to $(L^{i-d})_*$.
 \end{proof}

The niveau and the coniveau filtrations behave functorially with
respect to the action of correspondences:

 \begin{proposition} \label{respect} Let $\alpha \in \CH_{d+l}(X \times
   Y)$, then $\alpha_*N^jH_i(X) \subseteq N^{j+l}H_{i+2l}(Y)$ and
   $\alpha_*\N^jH_i(X) \subseteq \N^{j+l}H_{i+2l}(Y)$.
 \end{proposition}

 \begin{proof} The proof is straightforward from the definition of the
   niveau and coniveau filtrations.
\end{proof}

There is a non-degenerate bilinear form $Q_i$ on $H_i(X)$ given by
$$Q_i(\alpha,\beta) = \ \langle  L^{i-d}\alpha,\beta \rangle,$$ 
where $\langle  - , - \rangle \ : H_i(X) \otimes_\Q H_{2d-i}(X) \r \Q$
denotes the cup product.  The primitive part of $H_i(X)$ is defined to
be $H_i(X)_{\mathrm{prim}} := \ker \big(L^{i-d+1} : H_i(X) \rightarrow
H_{2d-i-2}(X)\big)$. We have the Lefschetz decomposition formula
$$H_i(X) = \bigoplus_{0 \leq 2r \leq 2d-i} L^r
H_{i+2r}(X)_{\mathrm{prim}}.$$ The Hodge index theorem states that the
above decomposition is orthogonal for $Q_i$ and that the sub-Hodge
structure $L^r H_{i+2r}(X,\Q)_{\mathrm{prim}}$ endowed with the form
$(-1)^rQ_i$ is polarized. We let $p_{i,r}$ denote the orthogonal
projector
$$H_i(X,\Q) \r L^rH_{i+2r}(X,\Q)_{\mathrm{prim}} \r H_i(X,\Q).$$
Polarized Hodge structures have the following well-known
particularity.

\begin{lemma} \label{subpol} Let $H$ be a polarized rational Hodge
  structure and let $K$ be a sub-Hodge structure. Then the pairing on
  $H$ remains non-degenerate after restriction to $K$.
\end{lemma}
\begin{proof} This is because the associated hermitian form on
  $H\otimes \C$ remains non-degenerate when restricted to the pieces
  $K_{p,q}$ of the Hodge decomposition of $K\otimes \C$.
\end{proof}

\begin{proposition} \label{pairing} If $X$ satisfies property B, then,
  for all $i$ and all $j$, the cup product pairings $${\widetilde{N}}^j
  H_i(X) \otimes {\widetilde{N}}^{d-i+j} H_{2d-i}(X) \r \Q \ \ \
  \mathrm{and} \ \ \ {{N}}^j H_i(X) \otimes {{N}}^{d-i+j} H_{2d-i}(X)
  \r \Q$$ are non-degenerate.
\end{proposition}
\begin{proof} The arguments that follow work equally well for the
  filtration $N$.

  First we have $L^{d-i}_* \widetilde{N}^j H_i(X) =
  \widetilde{N}^{d-i+j} H_{2d-i}(X)$. Indeed, if $i \leq d$, then, by
  Proposition \ref{respect}, $L^{d-i}_* \widetilde{N}^j H_i(X)
  \subseteq \widetilde{N}^{d-i+j} H_{2d-i}(X)$.  Because $L^{d-i}_*$
  is invertible as a correspondence, the reverse inclusion holds. It
  is thus enough to show that $(Q_i)|_{\widetilde{N}^jH_i(X)}$ is
  non-degenerate.

  Secondly, Kleiman showed \cite[Theorem 4.1.(3)]{Kleiman} that if $X$
  satisfies B, then the projectors $p_{i,r}$ are induced by algebraic
  correspondences $P_{i,r} \in \CH_d(X \times X)$.  Therefore,
  $(p_{i,r})_* \widetilde{N}^jH_i(X) = (P_{i,r})_*
  \widetilde{N}^jH_i(X)$ and, by Proposition \ref{respect}, we have
  $(p_{i,r})_* \widetilde{N}^jH_i(X) \subseteq \widetilde{N}^jH_i(X)$.
  Hence $$\widetilde{N}^jH_i(X) = \bigoplus_{0 \leq 2r \leq 2d-i}
  (p_{i,r})_* \widetilde{N}^jH_i(X).$$ Therefore, thanks to Lemma
  \ref{subpol}, the form $Q_i$ restricts to a non-degenerate form on
  $\widetilde{N}^jH_i(X)$.
\end{proof}

\begin{remark} \label{niveausplitting} The proof of the proposition
  shows in particular that if $X$ satisfies property B, then both
  filtrations $N$ and $\N$ on $H_i(X)$ split into orthogonal pieces
  for the form $Q_i$, e.g. as $H_i(X) = \bigoplus \Gr_{\N}^j H_i(X)$.
  Moreover we see that if $H_i(X)$ is polarized with respect to the
  choice of the ample line bundle $L$, then the splitting of the
  coniveau (resp. niveau) filtration $N$ (resp. $\N)$ on $H_i(X)$
  induced by the polarization identifies canonically with the
  splitting induced by the bilinear form $Q_i$.
\end{remark}

\subsection{Proof of Theorem \ref{Ksupport}} \label{Ksupportsection}

We start with two well-known lemmas.

\begin{lemma}[Lemma 5 in \cite{Voisin2}] \label{image} Let $H$ and
  $H'$ be rational Hodge structures endowed, respectively, with a
  polarization $Q$ and $Q'$.  Let $\gamma : H \r H'$ be a morphism of
  Hodge structures. Then $$\im (\gamma) = \im (\gamma \circ
  \gamma^\vee),$$ where $H$ is identified with its dual $H^\vee$ via
  the polarization $Q$ it carries.
\end{lemma}

\begin{proof}
  By linear algebra, we have $\dim_\Q H = \dim_\Q \ker \gamma +
  \dim_\Q \im \gamma^\vee$. Therefore, it is enough to prove that the
  subspace $\ker \gamma \cap \im \gamma^\vee$ is zero.  The subspace
  $\ker \gamma \cap \im \gamma^\vee$ actually defines a sub-Hodge
  structure of $H$ and we claim that the pairing $Q$ restricts to zero
  on this subspace. Indeed, let $z$ and $\gamma^\vee y$ be elements in
  $\ker \gamma \cap \im \gamma^\vee$. Then we have
   $$Q(\gamma^\vee y , z) \ = \ Q'(y, \gamma z) \ = \ Q'(y,0) \
   = 0.$$ By Lemma \ref{subpol}, $\ker \gamma \cap \im \gamma^\vee =
   \{0\}$.
 \end{proof}

 \begin{lemma}[Lemma 7 in \cite{Charles}] \label{lowdeg} Let $Z$ be a
   smooth projective variety of dimension $\leq k$ endowed with an
   ample line bundle $L_Z$. If $Z$ satisfies properties $\mathrm{B}_l$
   for all $l \leq 2\dim Z - k - 2$, then for all $r \geq 0$ there
   exists a correspondence whose action on $H_k(Z)$ is the orthogonal
   projector $H_k(Z) \r L_Z^r H_{k+2r}(Z)_{\mathrm{prim}} \r H_k(Z)$.
 \end{lemma}
 \begin{proof} By induction on $r$ it is enough to prove that there is
   a correspondence whose action on $H_k(Z)$ is the projector $H_k(Z)
   \r L_Z H_{k+2}(Z) \r H_k(Z)$. For this purpose, let $\alpha$ be a
   correspondence in $\CH_{k+2}(Z \times Z)$ inducing the inverse to
   the Lefschetz isomorphism $L_Z^{k +2 -\dim Z} : H_{k+2}(Z) \r
   H_{2\dim Z - k -2}(Z)$. Then $L_Z \circ \alpha \circ L_Z^{k+1-\dim
     Z}$ is the required correspondence because $
   H_{k}(Z)_{\mathrm{prim}}$ is the kernel of the action of
   $L_Z^{k+1-\dim Z}$ on $H_k(Z)$.
 \end{proof}

 We fix $i$ and $j$ with $j \geq 1$. By definition of the niveau
 filtration $\widetilde{N}$, let $Z$ be a smooth projective variety of
 pure dimension $i-2j$ and let $\Gamma \in \CH_{i-j}(Z\times X)$ be a
 correspondence such that $$ {\widetilde{N}}^j H_i(X) = \im \big(
 \Gamma_* : H_{i-2j}(Z) \r H_i(X)\big).$$ Because $X$ satisfies B,
 there is a cycle $\pi_i \in \CH_d(X\times X)$ whose homology class is
 the projector $H_*(X) \twoheadrightarrow H_i(X) \hookrightarrow
 H_*(X)$ (see \cite{Kleiman}).  Hence, considering $\pi_i \circ
 \Gamma$ instead of $\Gamma$ we can assume that $$ {\widetilde{N}}^j
 H_i(X ) = \im \big( \Gamma_* : H_{*}(Z ) \r H_*(X )\big).$$ Moreover,
 for weight reasons, $\Gamma$ acts as zero on $H_k(Z )$ for $k \neq
 i-2j$.

 Supposing that we have already exhibited the cycles $\pi_{i,k} \in
 \CH_d(X \times X)$ for $k>j$, we are going to construct a cycle
 $\pi_{i,j}$ as in the statement of the theorem. Note that, for $k
 >i/2$, $\pi_{i,k}=0$ is suitable. Consider then the cycle $\gamma :=
 \big(\pi_i - \sum_{k>j}\pi_{i,k} \big) \circ \Gamma$. The induced
 morphism $\gamma_* : H_*(Z) \r H_*(X)$ has $\Gr_{\N}^jH_i(X)$ for
 image, that is the orthogonal complement of ${\N}^{j+1}H_i(X)$ inside
 ${\N}^jH_i(X)$ for the form $Q_i$.

 Write $L$ for the Lefschetz isomorphism $L : H_i(X)
 \stackrel{\simeq}{\longrightarrow} H_{2d-i}(X)$, which is assumed to
 be induced by an algebraic correspondence. Because $\gamma_*H_*(Z) =
 \Gr_{\N}^jH_i(X)$ and by Remark \ref{niveausplitting}, we see that
 the transpose of $\gamma$ induces a morphism ${}^t\gamma_* \circ L :
 H_i(X) \r H_{i-2j}(Z)$ which is zero on the orthogonal complement of
 $\Gr_{\N}^jH_i(X)$ (for the form $Q_i$ or for the polarization
 induced by $L$; these split $\N$ the same way). Consider now the
 Lefschetz decomposition $H_{i-2j}(Z) = \bigoplus_{r\geq 0} L_Z^r
 H_{i-2j+2r}(Z)_{\mathrm{prim}}$ and write $p_r$ for the orthogonal
 projector onto $ L_Z^r H_{i-2j+2r}(Z)_{\mathrm{prim}}$. There are two
 cases at hand.

 If we are assuming that $N^{j+1}H_i(X) = \N^{j+1}H_i(X)$, then $\gamma
 \circ p_r =0$ unless $r=0$. Indeed, if $x \in L_Z^r
 H_{i-2j+2r}(Z)_{\mathrm{prim}}$, then $\gamma_*(x) = (\gamma \circ
 L_Z^r)_*(y)$ for some $y \in H_{i-2j+2r}(Z)_{\mathrm{prim}}=
 H^{i-2j-2r}(Z)_{\mathrm{prim}}$ so that $\gamma_*(x) \in
 N^{j+r}H_i(X)$. Thus, if $r>0$, then $\gamma_*(x) \in N^{j+1}H_i(X) =
 \N^{j+1}H_i(X)$ and hence $\gamma_*(x) = 0$. Therefore, the action of
 $\gamma$ on $H_*(Z)$ factors through the polarized Hodge structure
 $H_{i-2j}(Z)_{\mathrm{prim}}$. We then let $\varphi := \gamma \circ
 {}^t\gamma \circ L \in \CH_d(X \times X)$.

 If we are assuming that $Z$ satisfies property $\mathrm{B}_l$ for all
 $l \leq i-2j-2$, then the morphism $s_Z := \sum (-1)^rp_r$ on
 $H_{i-2j}(Z)$ is induced by a correspondence thanks to Lemma
 \ref{lowdeg} and $\langle -,s_Z(-)\rangle$ defines a polarization on
 $H_{i-2j}(Z)$. We then let $\varphi := \gamma \circ s_Z \circ
 {}^t\gamma \circ L \in \CH_d(X \times X)$.

 Either way, Lemma \ref{image} gives $$\varphi_*H_*(X) =
 \gamma_*H_*(Z) = \Gr_{\N}^jH_i(X).$$

 Therefore, $\varphi$ restricts to an isomorphism of the finite
 dimensional $\Q$-vector space $\Gr_{\widetilde{N}}^jH_i(X )$.  By the
 theorem of Cayley--Hamilton, there exists a polynomial $P \in \Q[X]$
 such that $(\varphi|_{\Gr_{\widetilde{N}}^jH_i(X)})^{-1} =
 P(\varphi|_{\Gr_{\widetilde{N}}^jH_i(X)})$. Then let $\psi$ be the
 correspondence $P(\varphi)$. The composite $\pi_{i,j} := \psi \circ
 \varphi$ thus induces the identity on the rational Hodge structure $
 \Gr_{\widetilde{N}}^jH_i(X )$ and is zero on its orthogonal
 complement, i.e. it is the required cycle.

 Finally, we set $\pi_{i,0} = \pi_i - \sum_{j \geq 1} \pi_{i,j}$ and
 this completes the proof of Theorem \ref{Ksupport}. \qed \medskip

 We now show that the cycles of Theorem \ref{Ksupport} coincide with
 the ones given by Kahn--Murre--Pedrini \cite[\S7.7]{KMP}.

 \begin{proposition}
   \label{Ksupportsection-prop} Let $X$ be a smooth projective variety
   that satisfies $(\star)$. Then the cycles $\pi_{i,j}$ of Theorem
   \ref{Ksupport} coincide modulo homological equivalence with the
   ones given in \cite[\S7.7]{KMP}
 \end{proposition}
 \begin{proof} The construction in \cite{KMP} is performed under the
   only assumption that $X$ satisfies B. Thus, to start with, suppose
   merely that $X$ satisfies B. Since we are working in characteristic
   zero, it follows \cite{Kleiman} that homological and numerical
   equivalence coincide on $X$ and $X \times X$. It follows from
   Jannsen's theorem \cite{Jannsen3} that the homological motive
   $h^{\mathrm{hom}}(X)$ and its Tate twists are contained in a full
   semisimple abelian subcategory of motives for homological
   equivalence. We have as in \cite[7.7.3]{KMP} a unique decomposition
 $$h^{\mathrm{hom}}_i (X) = \bigoplus_j h^{\mathrm{hom}}_{i,j}(X),$$
 where $h^{\mathrm{hom}}_{i,j}(X)(-j)$ is effective while
 $h^{\mathrm{hom}}_{i,j} (X)(-j-1)$ has no non-zero direct summand
 that is effective. Any direct summand $M$ of $h^{\mathrm{hom}}_i
 (X)$ for which $M(-j)$ is effective is then contained in $\bigoplus_{
   j' \geq j} h^{\mathrm{hom}}_{i,j'} (X)$. Applying the homology
 functor $H$ to the above decomposition gives a similar decomposition
 of $H_i(X)$ into subspaces $ H_{i,j}(X)$. The $H_{i,j}(X)$ are
 mutually orthogonal, because the polarization on $H_i(X)$ arises from
 a polarization on $h^{\mathrm{hom}}_i (X)$.  Since
 $h^{\mathrm{hom}}_{i,j} (X)(-j)$ is effective, there is, for some $Y$,
 a retraction $h^{\mathrm{hom}}(Y )(j) \r h^{\mathrm{hom}} _{i,j}
 (X)$.  If $\Gamma : h^{\mathrm{hom}}(Z)(j) \r h^{\mathrm{hom}}(X)$ is
 its composite with the embedding, then $\Gamma_* : H_{i-2j}(Z) \r
 H_i(X)$ has image $H_{i,j}(X)$. Thus, $$(1) \hspace{140pt}
 \bigoplus_{j'\geq j} H_{i,j'}(X) \subseteq \N^jH_i(X).  \hspace{140pt}
 {}$$

 Suppose now that $(\star)$ holds for $X$ and let $\pi_{i,j}$ be a
 cycle as in Theorem \ref{Ksupport} : its action on $H_*(X)$ induces
 the orthogonal projection on $\Gr_{\widetilde{N}}^jH_i(X)$ and there
 are a smooth projective variety $Z$ and correspondences $f \in
 \CH_{i-j}(Z \times X)$ and $g \in \CH_{d-j}(X \times Z)$ such that
 $\pi_{i,j} = f \circ g$. In particular, $f \circ g$ is an idempotent
 modulo homological equivalence. The image of $f \circ g$ is a
 homological motive $M$ which is a direct summand of
 $h^{\mathrm{hom}}(X)$.  The composite of the map
 $h^{\mathrm{hom}}(Z)(j) \r h^{\mathrm{hom}}(X)$ induced by $f$ with
 the projection $\h^\hom(X) \r M$ is a retraction, and thus $M$ is a
 direct summand of $h^{\mathrm{hom}}(Z)(j)$.  It follows that $M$ is
 contained in $\bigoplus_{ j' \geq j} h^{\mathrm{hom}}_{i,j'} (X)$.
 Applying $H$ then shows that the image $\Gr_{\widetilde{N}}^jH_i(X)$
 of $(f \circ g)_*$ is contained in the left-hand side of (1). Thus
 (1) is an equality. Both the $\Gr_{\widetilde{N}}^jH_i(X)$ and the
 $H_{i,j}(X)$ therefore give orthogonal splittings of the filtration
 $\N$ of $H_i(X)$, so that
 $$\Gr_{\widetilde{N}}^jH_i(X) = H_{i,j}(X)$$ as wanted.
\end{proof}

Since the decomposition of Kahn--Murre--Pedrini is unique modulo
homological equivalence, it follows that the idempotents $\pi_{i,j}$
when considered as endomorphisms of the homological motive $h^\hom(X)$
are unique. In particular, the cycles $\pi_{i,j}$ are central
endomorphisms of $h^\hom(X)$.

\begin{proposition} \label{independent-L} Let $X$ be a smooth
  projective variety that satisfies $(\star)$. Then the splitting of
  the niveau filtration on $H_i(X)$ is independent of the choice of
  polarization.
 \end{proposition}
 \begin{proof}
   Let's actually show that any two splittings of the niveau
   filtration on $H_i(X)$ that are defined by algebraic cycles
   coincide. This will prove the proposition because a splitting
   induced by the choice of a polarization is induced by algebraic
   cycles by Theorem \ref{Ksupport}. Let $\pi'_{i,j}$ be cycles in
   $\CH_d(X \times X)$ (not necessarily satisfying the conditions of
   Theorem \ref{Ksupport}) such that the endomorphism $(\pi'_{i,j})_*$
   of $H_i(X)$ induced by the cycles $\pi_{i,j}'$ are idempotents
   which give a splitting (not necessarily defined by a polarization)
   of the niveau filtration.  Then it is easily seen that the action
   of $\pi'_{i,j}$ on $H_i(X)$ coincides with the action of
   $\pi_{i,j}$, using descending induction on $j$ and the fact that
   two idempotent endomorphisms of a vector space coincide when they
   commute with one another and have the same image.
 \end{proof}

\subsection{Proof of Theorem \ref{action}}

\paragraph{A standard lifting lemma, proof of 4. in Theorem \ref{action}.}
\begin{lemma}[Particular case of \cite{Jannsen}, Lemma 3.1]
  \label{liftproj} Let $X$ be a variety satisfying $(\star \star)$.
  Then the following statements hold.
  \begin{enumerate}
  \item Let $p \in \CH_d(X \times X)$ be an idempotent. If
    $p_*H_*(X)=0$ then $p=0$.
  \item Let $c_1, \ldots, c_n \in \CH_d(X \times X)$ be correspondences
    such that $cl(c_i) \in H_*(X\times X) = \End(H_*(X))$ define
    mutually orthogonal idempotents adding to the identity. Then there
    exist mutually orthogonal idempotents $p_1, \ldots, p_n \in \CH_d(X
    \times X)$ adding to the identity in $\CH_d(X \times X)$ such that
    $cl(p_i) = cl(c_i)$ for all $1 \leq i \leq n$. Moreover, any two
    such choices $\{p_1, \ldots , p_n\}$ and $\{p'_1, \ldots , p'_n\}$
    are conjugate by an element lying above the identity, i.e. there
    exists a nilpotent correspondence $n \in \CH_d(X\times X)$ such
    that $p'_i = (1+n)\circ p_i \circ (1 + n)^{-1}$ for all $1 \leq i
    \leq n$.
    \end{enumerate}
\end{lemma}

Therefore if $X$ is a variety satisfying both $(\star)$ and $(\star
\star)$ then the correspondences $\pi_{i,j}$ of Theorem \ref{Ksupport}
can be chosen to be idempotents adding to $\Delta_X \in \CH_d(X \times
X)$ and any two such choices are conjugate. Let's fix such a choice
and write such idempotents $\Pi_{i,j}$.  As an immediate consequence
of the above lemma, we settle 4. of Theorem \ref{action} :

\begin{proposition} \label{vanish} If $X$ satisfies $(\star)$ and
  $(\star \star)$ and if $\Gr_{\widetilde{N}}^jH_{i}(X)=0$ for some
  integers $i$ and $j$, then $\Pi_{i,j}=0$.
\end{proposition}

\begin{proof}
  The identity $\Gr_{\widetilde{N}}^jH_{i}(X)=0$ means that the
  homology class of the idempotent $\Pi_{i,j}$ is zero. Part (1) of
  Lemma \ref{liftproj} implies that $\Pi_{i,j}=0$.
\end{proof}

\paragraph{Three lemmas.}

\begin{lemma} \label{proj} Let $X$ be a variety that satisfies $(\star
  \star)$.  Let $\pi$ be an idempotent and $\gamma$ be a
  correspondence in $\CH_d(X \times X)$, both acting on $\CH_*(X)$. Let
  also $S$ be a subset of $\ker \gamma$ which is stabilized by all
  nilpotent correspondences in $\CH_d(X \times X)$.  Then, if $\pi$ and
  $\gamma$ have same homology class, $S \subseteq \ker \pi$.
\end{lemma}

\begin{proof}
  By Lemma \ref{liftproj}, there exists a nilpotent correspondence $n
  \in \CH_d(X \times X)$ such that $\pi = \gamma + n$. Let $N$ be the
  nilpotent index of $n$.  Let $x\in S \subseteq \ker \gamma$. Then
  $\pi(x) = (\gamma + n)(x)$. Therefore, $\pi^{\circ N} (x) = (\gamma
  + n)^{\circ N} (x)$. By assumption $n$ stabilizes $S$ which is a
  subset of the kernel of $\gamma$.  Moreover, $n^{\circ N}=0$. Hence
  $\pi(x) = 0$, that is $x \in \ker \pi$.
\end{proof}

\begin{lemma} \label{proj2} Let $X$ be a variety that satisfies
  $(\star \star)$.  Let $\pi$ and $\pi'$ be two idempotents in $\CH_d(X
  \times X)$, both acting on $\CH_*(X)$. Suppose $\ker \pi'$ is
  stabilized by all nilpotent correspondences in $\CH_d(X \times X)$.
  Then, if $\pi$ and $\pi'$ have same homology class, $\ker \pi = \ker
  \pi'$.
\end{lemma}

\begin{proof}
  The previous lemma shows that $\ker \pi' \subseteq \ker \pi$. By
  Lemma \ref{liftproj}, there exists a nilpotent correspondence $n \in
  \CH_d(X \times X)$ such that $\pi = (1+n) \circ \pi' \circ
  (1+n)^{-1}$. Therefore, $\ker \pi = (1+n) \big( \ker \pi' \big)$ and
  we conclude that $\ker \pi \subseteq \ker \pi'$ since by assumption
  $n$ stabilizes $\ker \pi'$.
\end{proof}

\begin{lemma} \label{supportlemma} Let $X$ be a variety that satisfies
  $(\star \star)$ and let $\Pi$ be an idempotent in $\CH_d(X \times X)$
  which is homologically equivalent to a correspondence $\pi$ that
  factors through a smooth projective variety $Z$ : $\pi=f\circ g$ for
  $f \in \CH^{d-j}(Z \times X)$ and $g\in \CH_{d-j}(X \times Z)$.  Then
  there exist $f_1,\ldots,f_m \in \CH^{d-j}(Z \times X)$ and $g_1,
  \ldots, g_m\in \CH_{d-j}(X \times Z)$ such that $\Pi = \sum f_i \circ
  g_i$.
\end{lemma}
\begin{proof}
  There is a nilpotent correspondence $n \in \CH_d(X \times X)$ such
  that $\Pi= \pi +n$. Let $N$ be large enough so that $n^{\circ N}=0$.
  Then $\Pi = \Pi^{\circ N} = (\pi+n)^{\circ N}$. Expanding this last
  term, we see that each summand factors through $\pi$.
\end{proof}

\paragraph{Proof of 5. and 6. in Theorem 2}

The base field $k$ is an algebraically closed subfield of $\C$. In
\cite[Th. 3]{Vial3}, we proved the following (Proposition
\ref{pairing} shows that the conditions on $X$ in \textit{loc. cit.}
are met when $X$ satisfies property B). For the sake of completeness, we
provide a proof in the specific situation when $X$ satisfies B.

\begin{theorem} \label{Chowproj} Let $X$ be a smooth projective
  variety of dimension $d$ over $k$ that satisfies B. Then there
  exists a set of mutually orthogonal idempotents
  $\widetilde{\Pi}_{i,\lfloor i/2\rfloor} \in \CH_d(X \times X)$ $(0
  \leq i \leq 2d)$ that factor through curves and whose homology
  classes are the idempotents $H_*(X) \twoheadrightarrow
  \widetilde{N}^{\lfloor i/2\rfloor} H_{i}(X) \hookrightarrow H_*(X)$.
  Moreover if $\widetilde{\Pi}_{i,\lfloor i/2\rfloor}$ is any
  idempotent that factors through a curve and that induces the
  orthogonal projection $H_*(X) \twoheadrightarrow
  \widetilde{N}^{\lfloor i/2\rfloor} H_{i}(X) \hookrightarrow H_*(X)$,
  then we have
  \begin{center}
  $\CH_i(X)_{\hom} = \ker \big(\widetilde{\Pi}_{2i,i} : \CH_i(X) \r
  \CH_i(X)\big)$ and \medskip

 $\ker \big( AJ_i : \CH_i(X)_{\alg} \r J_i(X)\otimes \Q \big) =
  \ker \big(\widetilde{\Pi}_{2i+1,i} : \CH_i(X)_{\alg} \r
  \CH_i(X)_{\alg} \big).$ \end{center}
\end{theorem}
\begin{proof}
  First we show that there exist mutually orthogonal idempotents
  $\widetilde{\Pi}_{i,\lfloor i/2\rfloor}$ with the properties that
  the $\widetilde{\Pi}_{i,\lfloor i/2\rfloor}$ factor through a curve
  and that $(\widetilde{\Pi}_{i,\lfloor i/2\rfloor})_*H_*(X)=
  \widetilde{N}^{\lfloor i/2\rfloor} H_{i}(X)$. Note that for each $i$
  there is a $Z_i$ of dimension $1$ and a morphism $$q_i :
  h^{\mathrm{hom}}(Z_i)(\lfloor i/2 \rfloor) \r h^{\mathrm{hom}}(X)$$
  such that $(q_i)_*$ has image $ \N^{\lfloor i/2 \rfloor}H_i(X)$.
  Since $X$ and the $Z_i$ satisfy B, there is a full subcategory of
  the category of motives for homological equivalence containing
  $h^{\mathrm{hom}}(X)$ and the $h^{\mathrm{hom}}(Z_i)$ and their Tate
  twists which is semisimple abelian by Jannsen's theorem
  \cite{Jannsen3}, and on which $H$ is therefore exact.  Thus $q_i$
  has an image $\overline{M}_i$, and $H_*(\overline{M}_i) =
  \N^{\lfloor i/2 \rfloor}H_i(X)$. Since the polarization on $H_*(X)$
  arises from a polarization on $h^{\mathrm{hom}}(X)$, the
  decomposition of $H_*(X)$ as the sum of $ \N^{\lfloor i/2
    \rfloor}H_i(X)$ and its orthogonal complement arises from a
  similar direct sum decomposition of $h^{\mathrm{hom}}(X)$.  Thus,
  the embedding $\overline{s}_i : \overline{M}_i \r
  h^{\mathrm{hom}}(X)$ has a left inverse $\overline{r}_i :
  h^{\mathrm{hom}}(X) \r \overline{M}_i$ such that $(\overline{s}_i
  \circ \overline{r}_i)_*$ is the orthogonal projection onto $
  \N^{\lfloor i/2 \rfloor}H_i(X)$.  By semisimplicity, the image
  $\overline{M}_i$ of $q_i$ is a direct summand of
  $h^{\mathrm{hom}}(Z_i)(\lfloor i/2 \rfloor)$.  Since $h(Z_i)$ is
  finite-dimensional in the sense of Kimura, there exists a direct
  summand ${M}_i$ of $h(Z_i)(\lfloor i/2 \rfloor)$ lying above
  $\overline{M}_i$. The direct sum $M$ of the $M_i$ is then finite
  dimensional, and lies above the direct sum $\overline{M}$ of the
  $\overline{M}_i$. Write $\overline{s} : \overline{M} \r
  h^{\mathrm{hom}}(X)$ and $\overline{r} : h^{\mathrm{hom}}(X) \r
  \overline{M}$ for the morphisms with respective components the
  $\overline{s}_i$ and the $\overline{r}_i$. Then $\overline{r}$ is
  left inverse to $\overline{s}$.  If $s : M \r h(X)$ is a lifting of
  $\overline{s}$, then there is by finite-dimensionality of $M$ a
  lifting $r : h(X) \r M$ of $\overline{r}$ which is left inverse to
  $s$.  Write $s_i : M_i \r h(X)$ and $r_i : h(X) \r M_i$ for the
  respective components of $s$ and $r$. Then $r_i$ is left inverse to
  $s_i$, and we may take as $\widetilde{\Pi}_{i,\lfloor i/2\rfloor}$
  the idempotent $s_i\circ r_i$. That the $\widetilde{\Pi}_{i,\lfloor
    i/2\rfloor}$ are mutually orthogonal is clear, and
  $\widetilde{\Pi}_{i,\lfloor i/2\rfloor}$ factors through
  $h(Z_i)(\lfloor i/2 \rfloor)$ because it factors through the direct
  summand $M_i$ of $h(Z_i)(\lfloor i/2 \rfloor)$.

  The rest of the theorem then follows as in the proof of Lemma
  \ref{D} or \ref{B'} by functoriality of the cycle class maps and of
  the Abel--Jacobi maps with respect to the action of correspondences.
\end{proof}

If, moreover, $X$ satisfies $(\star \star)$, then it is possible to
prove that the above kernels do not depend on the choice of a lift
${\Pi}_{i,\lfloor i/2\rfloor}$ for the idempotent $\pi_{i,\lfloor
  i/2\rfloor}$.

\begin{proposition} \label{kerind} If $X$ satisfies B and $(\star
  \star)$ then
  $$\CH_i(X)_{\hom} = \ker \big({\Pi}_{2i,i} : \CH_i(X) \r
  \CH_i(X)\big)$$ for any choice of idempotent $\Pi_{2i,i}$
  such that $({\Pi}_{2i,i})_*H_{*}(X) =
  \widetilde{N}^iH_{2i}(X)$.
\end{proposition}

\begin{proof} By Theorem \ref{Chowproj}, $\CH_i(X)_{\hom} = \ker
  \big(\widetilde{\Pi}_{2i,i} : \CH_i(X) \r \CH_i(X)\big)$. The group
  $\CH_i(X)_{\hom}$ is stabilized by the action of correspondences in
  $\CH_d(X\times X)$. Therefore, by Lemma \ref{proj2}, $\ker
  \big(\Pi_{2i,i} \big) = \ker \big(\widetilde{\Pi}_{2i,i} \big)$,
  where both idempotents act on $\CH_i(X)$.
\end{proof}

\begin{proposition} \label{kerind2} If $X$ satisfies B and $(\star
  \star)$, then
  $$\ker \big( AJ_i : \CH_i(X)_{\alg} \r
  J_i(X)\otimes \Q \big) = \ker \big(\Pi_{2i+1,i} : \CH_i(X)_{\alg} \r
  \CH_i(X)_{\alg} \big),$$ where $\Pi_{2i+1,i}$ is any idempotent in
  $\CH_d(X\times X)$ such that $(\Pi_{2i+1,i})_* H_*(X)= \N^i
  H_{2i+1}(X)$.
\end{proposition}

\begin{proof} The idempotents $\widetilde{\Pi}_{2i+1,i}$ and
  $\Pi_{2i+1,i}$ are homologous.  By Lemma \ref{proj2} it is thus
  enough to see that $\ker \big( AJ_i : \CH_i(X)_{\alg} \r
  J_i(X)\otimes \Q \big) = \ker \big(\widetilde{\Pi}_{2i+1,i} \big)$
  is stabilized by the action of correspondences in $\CH_d(X\times X)$.
  This is indeed the case by functoriality of the Abel--Jacobi map with
  respect to the action of correspondences.
\end{proof}

\paragraph{proof of 1., 2. and 3.  in Theorem \ref{action}.}

The idempotent $\Pi_{i,j}$ is homologically equivalent to a
correspondence ${\pi}_{i,j}$ that factors through a variety $Z_{i,j}$
of dimension $i-2j$. The action of this correspondence $\pi_{i,j}$ on
$\CH_l(X)$ factors through $\CH_{l-j}(Z_{i,j})$, and hence $\ker
{\pi}_{i,j} = \CH_l(X)$ when $l < j$ or $l>i-j$ for dimension reasons.
Clearly, $\CH_l(X)$ is stabilized by the action of any correspondence
in $\CH_d(X\times X)$. Therefore Lemma \ref{proj} applies, giving $\ker
{\Pi}_{i,j} \supseteq \ker {\pi}_{i,j} = \CH_l(X)$ for $l < j$ or
$l>i-j$. Alternately one could have used Lemma \ref{supportlemma}
directly.

In the case $l=i-j$ and $i<2l$, the action of ${\pi}_{i,j}$ on
$\CH_l(X)$ factors through $\CH_{i-2j}(Z_{i,j})=\CH^0(Z_{i,j})$. In
particular, $\pi_{i,j}$ acts trivially on $\CH_l(X)_\hom$. Because
$i<2l$, $\pi_{i,j}$ acts trivially on $ H_{2l}(X)$ and hence on $\im
(\CH_l(X) \r H_{2l}(X))$.  Therefore, by functoriality of the cycle
class map, ${\pi}_{i,j}$ sends $\CH_l(X)$ to $\CH_l(X)_{\hom}$. Hence,
${\pi}_{i,j} \circ {\pi}_{i,j}$ acts as zero on $\CH_l(X)$. Besides,
${\pi}_{i,j} \circ {\pi}_{i,j}$ is homologically equivalent to
${\Pi}_{i,j} \circ {\Pi}_{i,j} = {\Pi}_{i,j}$. By Lemma \ref{proj}, we
get $\ker {\Pi}_{i,j} \supseteq \ker \big( {\pi}_{i,j} \circ
{\pi}_{i,j} \big) = \CH_l(X)$.

Similarly, in the case $l+1=i-j$ and $i \leq 2l$, the action of
${\pi}_{i,j}$ on $\CH_l(X)$ factors through $\CH_{i-2j-1}(Z_{i,j})=
\CH^1(Z_{i,j})$. We have an isomorphism $AJ^1 : \CH^1(Z_{i,j})_\hom \r
J^1(Z_{i,j})$. Therefore, by functoriality of the Abel-Jacobi map on
homologically trivial cycles, ${\pi}_{i,j}$ acts trivially on $\ker
(\CH_l(X)_\hom \r J_l(X) )$. Now, by assumption on $i$, $j$ and $l$,
${\pi}_{i,j}$ acts trivially on $\im (\CH_l(X) \r H_{2l}(X))$ and on
$H_{2l+1}(X)$ and hence on $J_l(X)$. As such, ${\pi}_{i,j}$ sends
$\CH_l(X)$ to $\ker (\CH_l(X)_\hom \r J_l(X) )$. As before, we get $\ker
{\Pi}_{i,j} \supseteq \ker \big( {\pi}_{i,j} \circ {\pi}_{i,j} \big) =
\CH_l(X)$.  \qed \medskip

If $X$ satisfies $(\star)$ and $(\star \star)$, then a star in the diagram
below with coordinates $(i,j)$ $(0 \leq i \leq j)$ indicates that
$\mathrm{Gr}_{\widetilde{N}}^i H_{i+j}(X)$ does not induce a
``motivic'' action on $\CH_l(X)$, i.e. that $\Pi_{i+j,i}$ acts as zero
on $\CH_l(X)$. A bullet with coordinates $(i,j)$ $(0 \leq i \leq j)$
indicates that $\Pi_{i+j,i}$ should act as zero on $\CH_l(X)$ if one
believes in the BBM conjectures, see Proposition \ref{murre}.

$$\begin{array}{c | cc
    c c c cc c c c c c c} & 0 & 1 & & & & l & & & & & & & d\\ \hline
  0  &* & * & *  & * &* &*  &*  & \bullet  & \bullet  & \bullet   & \bullet \\
  1  &* & * & *  & * &* &* &*  & \bullet  & \bullet   & \bullet  \\
  &   *& *  & *  & * &* &*  &*  & \bullet  & \bullet  \\
  & * & * & *  & * &* &*  &*   & \bullet \\
  & * & * & *  & * &* &*  &*   & \\
  l &  * & *  &  *   &  *  &    *             &    \\
  & * & * & *  &  * &  *  &  & *  & *  & *  & *  & *  & * & *\\
  & \bullet & \bullet & \bullet & \bullet & & & * & * & * & *
  & *  & * & *\\
  &\bullet & \bullet  &  \bullet &  &    &  & *  & *  & *  & *  & *  & * & *\\
  & \bullet &  \bullet   &   &  &  &  & *  & *  & *  & *  & *  & * & *\\
  & \bullet  &   &   & &   &  & *  & *  & *  & *  & *  & * & *\\
  & &   &   &  &  &  & *  & *  & *  & *  & *  & * & *\\
  d    & &   &  & &    &  & *  & *  & *  & *  & *  & * & *\\
\end{array}$$   \medskip

We have ``symmetrized'' the diagram by including points with
coordinates $(i,j)$ with $i>j\geq 0$. As such, a point with
coordinates $(i,j)$ represents $\mathrm{Gr}_{\widetilde{N}}^{\min
  (i,j)} H_{i+j}(X)$. This way the diagram looks like a Hodge diamond.
The reason for doing so is the following.  Grothendieck's generalized
Hodge conjecture (GHC for short) predicts that the coniveau filtration
coincides with the Hodge coniveau filtration. Since the Lefschetz
standard conjecture is a particular instance of the Hodge conjecture,
Grothendieck's GHC for all varieties implies actually that the niveau
filtration should coincide with the Hodge coniveau filtration; see
Proposition \ref{filtration-coincide}. Therefore, if $H_{i+j}(X)$ is
given with a polarization then, under GHC, ${\widetilde{N}}^{\min
  (i,j)} H_{i+j}(X)$ is the sub-Hodge structure of $H_{i+j}(X)$
generated by the orthogonal complement inside $H_{i,j}(X) \oplus
H_{j,i}(X)$ of the intersection of $H_{i,j}(X) \oplus H_{j,i}(X)$ with
the complexification of the sub-Hodge structure of $H_{i+j}(X)$ spanned
by $H_{i+j,0}(X) \oplus \ldots \oplus H_{\max (i,j)+1, \min (i,j)}(X)
\oplus H_{\min (i,j), \max (i,j)+1}(X)\oplus \ldots \oplus H_{0,
  i+j}(X)$. Therefore, $\mathrm{Gr}_{\widetilde{N}}^{\min (i,j)}
H_{i+j}(X) = 0$ if and only if $H_{i,j}(X) \oplus H_{j,i}(X)$ is
included in the complexification of the sub-Hodge structure of $H_i(X)$
spanned by $H_{i+j,0}(X) \oplus \ldots \oplus H_{\max (i,j)+1, \min
  (i,j)}(X) \oplus H_{\min (i,j), \max (i,j)+1}(X)\oplus \ldots \oplus
H_{0, i+j}(X)$. In particular, if the Hodge number $h_{i,j}(X)$
vanishes, then $\mathrm{Gr}_{\widetilde{N}}^{\min (i,j)} H_{i+j}(X) =
0$.

\section{Some applications} \label{Some_applications}

In \S\S 2.1 and 2.2 we give direct applications of Theorem
\ref{action} and show explicitly how the niveau filtration on $H_*(X)$
reflects on the Chow groups $\CH_*(X)$ for those varieties $X$ that
satisfy $(\star)$ and $(\star \star)$. The most general result in this
direction is Proposition \ref{gen}. In particular, we relate (in the
spirit of Bloch's conjecture for surfaces) the support of the Chow
group of $0$-cycles of a variety $X$ satisfying $(\star)$ and $(\star
\star)$ to the support of its cohomology; see Proposition
\ref{support}. In \S 2.2, we give partial answers to questions raised
by Schoen and by Esnault and Levine. In \S \ref{Fermatexample} we are
interested in exhibiting a new example of a variety $Y$ for which we
can give a description of its Chow groups. This is achieved by
considering a finite quotient of a Fermat $4$-fold of degree $7$; see
Theorem \ref{example}.

\subsection{On the support of the Chow groups}

\begin{definition}
  The niveau number $g_{i,j}(X)$ of $X$ is equal to $\dim_\Q
  \Gr_{\widetilde{N}}^iH_{i+j}(X)$.
\end{definition}

In particular if the niveau numbers $g_{0,i}(X), g_{1,i-1}(X), \ldots,
g_{k,i-k}(X)$ vanish then $H_i(X) = \widetilde{N}^{k+1} H_i(X)$.
Consequently the Hodge numbers $h_{0,i}(X), \ldots, h_{k,i-k}(X)$
vanish. Here $h_{p,q}(X) := \dim_\C H^{d-p,d-q}(X)$.

\paragraph{$0$-cycles.}
Bloch and Srinivas \cite{BS} proved using a decomposition of the
diagonal that if $\CH_0(X \times \C)$ is supported in dimension $i$
(i.e. if there exists a closed subscheme $Z$ of $X \times \C$ of
dimension $i$ such that the map $\CH_0(Z) \r \CH_0(X \times \C)$ induced
by the inclusion of $Z$ inside $X \times \C$ is surjective), then
$H_k(X) = N^1H_k(X)$ for all $k >i$. This generalized the result of
Roitman \cite{Roitman} who considered the case $i=1$. In this latter
case we actually have a more precise result involving the niveau
filtration instead of the coniveau filtration.

\begin{proposition} \label{relevance} Let $X$ be a smooth projective
  complex variety. If $\CH_0(X)$ is supported on a curve, then $H_k(X)
  = \N^1H_k(X)$ for all $k >1$.
\end{proposition}
\begin{proof}
  By Bloch and Srinivas \cite{BS}, there exist a divisor $D$ on $X$
  and a curve $C$ in $X$ such that the diagonal $\Delta_X$ decomposes
  as $\Delta_X = \Gamma_1 + \Gamma_2 \in \CH_d(X \times X)$ with
  $\Gamma_1$ supported on $D \times X$ and $\Gamma_2$ supported on $X
  \times C$. Let $\widetilde{C}$ and $\widetilde{D}$ be
  desingularizations of $C$ and $D$ respectively. Then the action of
  $\Gamma_1$ on $H_k(X)$ factors through $H_{k-2}(\widetilde{D})$ and
  the action of $\Gamma_2$ on $H_k(X)$ factors through
  $H_{k}(\widetilde{C})$. This makes it possible to conclude that
  $H_k(X) = \N^1H_k(X)$ for all $k >1$.
\end{proof}

We now give a converse statement for those varieties that satisfy
$(\star)$ and $(\star$$\star)$. This gives a proof of Jannsen's
\cite[conjecture 3.3]{Jannsen} for varieties that satisfy $(\star)$
and $(\star$$\star)$, which is a slight improvement since Jannsen
prove gave a proof under the validity of the standard conjectures and
the BBM conjectures (the standard conjectures imply $(\star)$ and
Jannsen \cite{Jannsen} proved that the BBM conjectures imply that
property $(\star$$\star)$ holds for all smooth projective varieties).

\begin{proposition} \label{support} Let $X$ be a smooth projective
  variety satisfying $(\star)$ and $(\star \star)$. If $H_k(X)=
  \widetilde{N}^1H_k(X)$ for all $k > i$, then $\CH_0(X)$ is supported
  in dimension $i$.
\end{proposition}

\begin{proof}
  By Theorem \ref{action}, we have $$(\Delta_X)_* \CH_0(X) =
  ({\Pi}_{0,0} + {\Pi}_{1,0} + \ldots + {\Pi}_{i,0})_* \CH_0(X).$$ By
  Theorem \ref{Ksupport}, the idempotent ${\Pi}_{0,0} + {\Pi}_{1,0} +
  \ldots + {\Pi}_{i,0}$ is homologically equivalent to a cycle that
  factors through an $i$-dimensional variety. The proposition then
  follows from Lemma \ref{supportlemma}.
\end{proof}

\begin{remark}
  A surface $X$ satisfies $(\star)$. Moreover, if $X$ has vanishing
  geometric genus $p_g := h_{2,0}$, then, by the Lefschetz
  $(1,1)$-theorem, $H_2(X)= \widetilde{N}^1 H_2(X)$.  Therefore, we
  recover the following result due to Kimura \cite[Corollary
  7.7]{Kimura} : if $X$ is a surface with $p_g=0$ whose Chow motive is
  finite-dimensional, then the Bloch conjecture is true for $X$.
\end{remark}

\begin{remark}
  The generalized Bloch conjecture states that if the Hodge numbers
  \linebreak $h_{i+1,0}, h_{i+2,0}, \ldots, h_{d,0}$ of $X$ vanish,
  then $\CH_0(X)$ is supported in dimension $i$. Proposition
  \ref{support} shows that if Grothendieck's generalized Hodge
  conjecture holds then the generalized Bloch conjecture holds for
  those $X$ that satisfy $(\star \star)$ (the property $(\star)$ is
  automatically satisfied if Grothendieck's generalized Hodge
  conjecture holds).
\end{remark}

In the cases $i=0$ and $i=1$, the integral Chow group $\CH^\Z_0(X)$ of
$0$-cycles of $X$ can be computed explicitly.

\begin{proposition} \label{integralsupport} Let $X$ be a smooth
  projective variety over an algebraically closed field $k \subseteq
  \C$ satisfying $(\star)$ and $(\star \star)$.  If $H_i(X) =
  \widetilde{N}^1H_i(X)$ for all $i>1$, then $\CH^\Z_0(X) = \Z \oplus
  \mathrm{Alb}_X(k).$
\end{proposition}
\begin{proof} This is a refinement of a result of Voisin \cite[Theorem
  3]{Voisin2}. By Theorem \ref{action}, we have $\CH_0(X) = (\Pi_{0,0}
  + \Pi_{1,0})_* \CH_0(X).$ Still by Theorem \ref{action}, we have
  $\ker (\mathrm{alb}_X : \CH_0(X)_\hom \r \mathrm{Alb}_X(k)\otimes \Q)
  = \ker ((\Pi_{1,0})_* : \CH_0(X)_\hom \r \CH_0(X)_\hom)$. Thus $\ker
  (\mathrm{alb}_X) = 0$.  By Roitman's theorem \cite{Roitmantors}, the
  integral albanese map $\mathrm{alb}^\Z_X : \CH^\Z_0(X)_\hom \r
  \mathrm{Alb}_X(k)$ is an isomorphism on the torsion. It is also
  surjective, and hence the proposition.
\end{proof}

\paragraph{$l$-cycles.}

Lewis \cite{LewisMumford} and Schoen \cite{Schoen} proved
independently the following theorem.

\begin{theorem}[Lewis, Schoen] \label{LS} If $\Gr_{\N}^l H_{2l+k}(X)
  \neq 0$ for some $k\geq 2$ then $\CH_l(X)_\alg$ is not representable,
  i.e. there is no curve $C$ and no correspondence $\Gamma \in
  \CH_{l+1}(C \times X)$ such that $\Gamma_*\CH_0(C)_\alg =
  \CH_l(X)_\alg$.
\end{theorem}

\noindent The statement below, when considered in the case $p=1$ and
$q=0$, gives a partial converse for varieties $X$ satisfying $(\star)$
and $(\star \star)$.

\begin{proposition}[Generalization of Proposition \ref{support}]
  \label{gen} Let $X$ be a smooth projective variety satisfying
  $(\star)$ and $(\star \star)$ and let $l$ be a non-negative integer.
  If there exist non-negative integers $p$ and $q$ such that
  $g_{i,j}(X)=0$ when \vspace{1pt}

  (i) $i+j \leq 2l$ and $j>l+1$ \vspace{1pt}

  \noindent and when \vspace{1pt}

  (ii) $2l+1 \leq i+j \leq 2l+p$ and $i<l-q$ \vspace{1pt}

  \noindent and when \vspace{1pt}

  (iii) $ i+j > 2l+p$ and $i \leq l$ \vspace{3pt}

  \noindent then there exist a smooth projective variety $Z$ of
  dimension $p+2q$, a correspondence $\Gamma' \in \CH_{d+q-l}(X \times
  Z)$ and a correspondence $\Gamma \in \CH_{p+q+l}(Z \times X)$ such
  that $\Gamma \circ \Gamma'$ acts as the identity on $\CH_l(X)$. In
  particular $\CH_l(X) = \Gamma_*\CH_q(Z)$.
\end{proposition}

\begin{proof} Let $\Pi := {\Pi}_{2l,l} + ({\Pi}_{2l+1,l} + \ldots +
  {\Pi}_{2l+p,l}) + ({\Pi}_{2l+1,l-1} + \ldots + {\Pi}_{2l+p,l-1}) +
  \ldots + ({\Pi}_{2l+1,l-q} + \ldots + {\Pi}_{2l+p,l-q}).$ Under the
  assumptions on the niveau numbers, Theorem \ref{action} says
  $\CH_l(X) = \Pi_* \CH_l(X)$.  Theorem \ref{Ksupport} says that the
  idempotent $\Pi_{i,j}$ is homologically equivalent to a cycle of the
  form $f_{i,j} \circ g_{i,j}$ for some smooth projective variety
  $Z_{i,j}$ of dimension $i-2j$, some $f_{i,j} \in \CH_{i-j}(Z_{i,j}
  \times X)$ and some $g_{i,j} \in \CH_{d-j}(X \times Z_{i,j})$. The
  idempotents $\Pi_{i,j}$ appearing in the above sum all satisfy $i-2j
  \leq p+2q$ and $l-j \leq q$. Therefore, up to replacing each
  $Z_{i,j}$ by $Z_{i,j} \times \P^{p+2q-i+2j}$, Lemma
  \ref{supportlemma} implies the desired result.
\end{proof}

\subsection{Injectivity of some cycle class maps}

\begin{proposition} \label{aaa} Let $X$ be a smooth projective variety
  satisfying $(\star)$ and $(\star \star)$.  If the niveau numbers
  $g_{i,j}(X)$ vanish for \vspace{1pt}

  (i)  $j > l+1$ when $i + j \leq 2l$ \vspace{1pt}

  \noindent and for \vspace{1pt}

  (ii) $i \leq l$ when $i + j > 2l$, \vspace{3pt}

\noindent then $\CH_l(X)_{\hom} = 0$, i.e. the cycle class map $\CH_l(X)
\r H_{2l}(X)$ is injective.
\end{proposition}

\begin{proof} By assumption on the niveau numbers, Theorem
  \ref{action} shows that the only idempotent $\Pi_{i,j}$ acting non
  trivially on $\CH_l(X)$ is $\Pi_{2l,l}$. Therefore $\CH_l(X) =
  (\Pi_{2l,l})_* \CH_l(X)$.  Hence $\CH_l(X)_{\hom} = \ker
  ({\Pi}_{2l,l}) = 0$.
\end{proof}

\begin{remark} The $4$-fold $Y$ that will be considered in \S
  \ref{Fermatexample} satisfies the assumptions of the proposition
  with $l=2$.
\end{remark}

\begin{corollary} Let $X$ be a smooth projective variety satisfying
  $(\star)$ and $(\star \star)$. If the niveau numbers $g_{p,q}(X)$
  vanish for $p\neq q$ and $p \leq k$, then $\CH_l(X)_{\hom} = 0$ for
  all $l \leq k$ and all $l \geq d-k-2$.
\end{corollary}

\begin{proof} Note that the condition (i) of Proposition \ref{aaa}
  implies $i \leq l-2$.  It is thus easy to check that the vanishing
  condition on the niveau numbers implies, for all $l \leq k$, the
  vanishing conditions of Proposition \ref{aaa}. Therefore
  $\CH_l(X)_{\hom} = 0$ for all $l \leq k$. A generalized decomposition
  of the diagonal as done by Laterveer \cite{Laterveer} then shows
  that $\CH_l(X)_{\hom} = 0$ for all $l \geq d-k-2$.
\end{proof}

\begin{proposition} [Answers partially question 0.6 of \cite{Schoen}]
  \label{del} Let $X$ be a smooth projective variety satisfying
  $(\star)$ and $(\star \star)$. Suppose that there is an integer $l \geq
  0$ such that the niveau numbers $g_{i,j}(X)$ vanish if \vspace{1pt}

(i) $i+j \leq 2l$ and $ j > l+1$ \vspace{1pt}

\noindent and if \vspace{1pt}

(ii) $i+j = 2l+1$ and $ |i-j| > 1$ \vspace{1pt}

\noindent and if \vspace{1pt}

(iii) $i+j \geq 2l+2$ and $i \leq l$. \vspace{3pt}

\noindent Then $\CH_l(X)_{\hom } = \CH_l(X)_{\alg }$ and the Abel-Jacobi
map $AJ_l : \CH_l(X)_{\hom } \r J_l(X) \otimes \Q$ is injective.
\end{proposition}

\begin{proof}
  The assumptions made on the niveau numbers of $X$ and Theorem
  \ref{action} (points 1, 2, 3 and 4) give $\CH_l(X) =
  \big({\Pi}_{2l,l}+ {\Pi}_{2l+1,l}\big)_*\CH_l(X)$. The proposition
  then follows from points 5 and 6 of Theorem \ref{action}.
\end{proof}

\begin{corollary} [see question 1 of \cite{EL}]
  \label{del2} Let $X$ be a smooth projective variety satisfying
  $(\star)$ and $(\star \star)$. Suppose that there is an integer $s
  \geq 0$ such that the niveau numbers $g_{i,j}(X)$ vanish if

  (i) $i+j \leq 2s+2$ and $|i-j| >1$

  \noindent and if

  (ii) $i+j > 2s+2$ and $ i \leq  s$

  \noindent Then $\CH_l(X)_{\hom } = \CH_l(X)_{\alg }$ and the
  Abel--Jacobi maps $AJ_l : \CH_l(X)_{\hom } \r J_l(X) \otimes \Q$ are
  injective for $l = 0, \ldots, s$ and $l = d-s-1, \ldots, d$.
\end{corollary}

\begin{proof} It is not too difficult to check that the vanishing
  condition on the niveau numbers implies for all $l \leq s$ the
  vanishing conditions of Proposition \ref{del}. Therefore
  $\CH_l(X)_{\hom } = \CH_l(X)_{\alg }$ and the Abel--Jacobi maps $AJ_l :
  \CH_l(X)_{\hom } \r J_l(X) \otimes \Q$ are injective for all $l \leq
  s$. As before, a generalized decomposition of the diagonal as done
  by Laterveer \cite{Laterveer} makes it possible to conclude for $l
  \geq d-s-2$.
\end{proof}

The diagram below illustrates the assumptions made on the niveau
numbers of $X$ in the above corollary.  A star with coordinates
$(i,j)$ indicates that the niveau number $g_{\min (i,j),\max
  (i,j)}(X)$ is possibly non-zero. The absence of a star at the point
of coordinate $(i,j)$ indicates that $g_{\min(i,j),\max (i,j)}(X)=0$.

$$\begin{array}{c | c c @{ \ \,}c c@{\hskip-1mm } c @{}c c
    @{\hskip-3mm }c@{\hskip-3mm } c c @{}c@{} c}
  & 0 & 1 &  2  &        & s+1   &   &       & d-s-1  &   &   &
  \hskip-2mm
  d-1 & \ d\\ \hline
  0 & * & *   \\
  1 & * & *  & *   \\
  2 &    & *  & * & \ddots    \\
  &     &    & \ddots &   \ddots  &      *  &    \\
  s+1 &    &     &              &               *  & * & *  &  *  & *  \\
  & &   &              &                  & * & *  &       *       & *   \\
  &    &    &                           &     & * &  *&   *    &   *  \\
  d-s -1&  &   &    &    & *     & *   & * &  *  &  *           \\
  &    &    &   &     &        &     &   &     *  & \ddots  & \ddots  \\
  &     &   &      &     &   &     &   &      & \ddots &   *     &  * \\
  d-1   &    &   &     &   &   &     &   &     &     &   *   &  * & * \\
  d    &        &  &  &     &     &    &  &      &    &  &  *  & *
   \end{array}$$   \\

\begin{remark} In the case $s=d$, the condition $(\star)$ is
  automatically satisfied since it is always true that $N^{\lfloor i/2
    \rfloor +1}H_i(X) = \N^{\lfloor i/2 \rfloor+1}H_i(X) = 0$, cf.
  also \cite[Th. 4]{Vial3}.
\end{remark}

\subsection{A $4$-fold $Y$ of general type with $\CH^\Z_0(Y)=\Z$ }
\label{Fermatexample}

Koll\'ar showed \cite{Kollar} that any Fano variety $X$ satisfies
$\CH_0(X) = \Q$. A natural question is to ask whether or not there
exist varieties $X$ of general type with $\CH_0(X) = \Q$. This question
has a positive answer in the case of surfaces (see, for instance,
\cite{Barlow}, \cite{IM}, \cite{Godeaux}) and in the case of
threefolds \cite{Peters}. Given any two varieties $X$ and $Y$ which
satisfy $\CH_0(X) = \Q$ and $\CH_0(Y)=\Q$, it is easy to see that
$\CH_0(X\times Y) = \Q$ as well. Consider indeed $P$ (resp. $Q$) a
closed point of $X$ (resp. $Y$). Every $0$-cycle on $X$ is rationally
equivalent to a rational multiple of the class $[P]$ of $P$. It
follows that every $0$-cycle on $X \times Y$ is rationally equivalent
to a $0$-cycle of the form $[P] \times \alpha$ for some $0$-cycle
$\alpha \in \CH_0(Y)$. But then, $\alpha$ is rationally equivalent to a
rational multiple of $[Q]$. Every $0$-cycle on $X \times Y$ is
therefore rationally equivalent to a rational multiple of $[P \times
Q]$. If $X$ and $Y$ are of general type, then $X \times Y$ is of
general type. It is thus possible to exhibit, for all $d \geq 2$,
examples of varieties $X$ of dimension $d$ which are of general type
with $\CH_0(X)=\Q$.

More can be said. By Guletskii and Pedrini \cite{GP}, the Chow motive
of a complex surface $S$ with $\CH_0(S) = \Q$ is a sum of Lefschetz
motives. Therefore, the motive of the product of such surfaces is also a
sum of Lefschetz motives. Hence, the Chow groups of the product of
surfaces with $\CH_0(S) = \Q$ are finite-dimensional $\Q$-vector
spaces. If the surfaces involved in the product are all surfaces of
general type then the product is itself of general type.  It is thus
possible to exhibit, in any even dimension, an example of a variety
$X$ of general type with finite-dimensional Chow groups. Concerning
complex threefolds $W$ with $\CH_0(W)= \Q$, by \cite[Theorem 4 \& Rk
3.8]{Vial3} the motive of $W$ decomposes as a direct sum of Lefschetz
motives with the direct summand of a motive of curve tensored with the
Lefschetz motive. \medskip

Here we exhibit a new example of a complex fourfold $Y$ of general
type with $\CH_0(Y)=\Q$. We will also show that $AJ_1 : \CH_1(Y)_\alg \r
J_1(Y) \otimes \Q$ is not injective. Thus our variety $Y$ is not
isomorphic to the product of two surfaces of general type with
$\CH_0(S) = \Q$.  (Actually, it is simpler to note that the $H^4$ of a
product of two such surfaces is concentrated in Hodge bidegrees
$(2,2)$, but $H^{3,1}(Y) \neq 0$.) It is also not isomorphic to the
blow-up of the product of two surfaces of general type with $\CH_0(S) =
\Q$ along a smooth surface since such a variety has Picard number $>1$
whereas $Y$ will be seen to have Picard number $1$.  Let's mention
that there exist examples of varieties not of general type satisfying
the above : Schoen \cite[Theorem 0.5]{Schoen} showed that
hypersurfaces $X$ of dimension $d \geq 4$ and $d/2+1 \leq \deg X <
d+2$ have non injective $AJ_1 : \CH_1(X)_\alg \r J_1(X) \otimes
\Q$.\medskip

Let $X_m^n$ be the complex Fermat $n$-fold of degree $m$, that is $$
X_m^n= \{x_0^m + x_1^m + \ldots +x_{n+1}^m =0\} \subset \P_\C^{n+1}.$$
Such a variety satisfies property B (this is the case for any
hypersurface) and is dominated by a product of curves \cite{KS}. When
$m$ is a prime number, Ran \cite{Ran} and Shioda \cite{Shioda} proved
that the Hodge conjecture holds for $X_m^n$.

\noindent From now on $X$ denotes the complex Fermat fourfold of
degree $7$, i.e. $X=X_7^4 \subset \P_\C^{5}.$

\begin{theorem} \label{GHC} The Fermat fourfold $X$ satisfies
  Grothendieck's generalized Hodge conjecture. More is true :
  $\N^1H^4(X)$ is the largest sub-Hodge structure of $H^4(X)$ whose
  complexification is included in $H^{1,3}(X) \oplus H^{2,2}(X) \oplus
  H^{3,1}(X)$.
\end{theorem}

\begin{proof}
  We use Ran's and Shioda's technique (which is summed up in \cite[\S
  3 (13)]{Shioda2}) to prove
$$\N^1H^4(X) \cap H^4_{prim}(X,\Q)= (H^{1,3}(X)
\oplus H^{2,2}(X) \oplus H^{3,1}(X)) \cap H^4_{prim}(X,\Q).$$ From
 \cite[Theorem I]{Shioda}, we have a decomposition $$ \big( (H^{1,3}
\oplus H^{2,2} \oplus H^{3,1}) \cap H^4_{prim}(X,\Q) \big) \otimes_\Q
\C = \bigoplus_{\alpha \in \mathfrak{B}_7^4} V(\alpha),$$ where, if
$<a>$ denotes the least positive residue of $a$ modulo $7$,
$$\mathfrak{B}_7^4 = \{(a_0,a_1,a_2,a_3,a_4,a_5) \ | \ 1 \leq a_i \leq
6, \ \sum_{i=0}^5 <ta_i>= 14, \ 21 \ \mathrm{or} \ 28 \ \mathrm{for} \
\mathrm{all} \ t \in (\Z/7)^\times \}$$ and $V(\alpha)$ is a certain
subspace of $H^4_{prim}(X,\C)$ of dimension $1$. An element $\alpha =
(a_i) \in \mathfrak{B}_7^4$ is called \textit{decomposable} if
$a_i+a_j = 0$ mod $7$ for some $i \neq j$ and
\textit{quasi-decomposable} if $(a_0,a_1, a_2, a_3, <a_4 + a_5>)$
belongs to $\mathfrak{B}_7^3 = \{(b_0,b_1,b_2,b_3,b_4) \ | \ 1 \leq
b_i \leq 6, \ \sum_{i=0}^4 <tb_i>= 14 \ \mathrm{or} \ 21 \
\mathrm{for} \ \mathrm{all} \ t \in (\Z/7)^\times \}$ after a
permutation of the digits $a_i$. By Theorem II of \cite{Shioda}, there
is an isomorphism induced by a correspondence
$$[H^3(X^3_7) \otimes H^1(X^1_7)]^{\mu_7} \oplus
[H^2_{prim}(X^2_7) \otimes H^0(X^0_7)](-1)
\stackrel{\simeq}{\longrightarrow} H^4_{prim}(X^4_7).$$ If $\alpha \in
\mathfrak{B}_7^4$ is decomposable (resp. quasi-decomposable), then
$V(\alpha)$ corresponds to some $V(\beta) \otimes V(\gamma) \in
H^2_{prim}(X^2_7) \otimes H^0_{prim}(X^0_7)$ (resp. to some $V(\beta)
\otimes V(\gamma) \in H^3_{prim}(X^3_7) \otimes H^1_{prim}(X^1_7)$
with $\beta \in \mathfrak{B}_7^3)$. Therefore, if $\alpha$ is
decomposable, then it comes from a class supported on a surface, i.e.
$V(\alpha) \subset \N^1H^4(X,\C)$. If $\beta \in \mathfrak{B}_7^3$,
then, by Shioda \cite[\S 3 (13)]{Shioda2}, $V(\beta) \subset
H^3(X^3_7) $ is supported in codimension one.  Having in mind that
$N^1H_3(Z) = \N^1H_3(Z)$ for a $3$-fold $Z$, we get that $V(\beta)$
comes from a curve.  Therefore, if $\alpha$ is quasi-decomposable,
then $V(\alpha)$ comes from a class supported on a surface. Now, it
can be checked that any element of $\mathfrak{B}_7^4$ is either
decomposable or quasi-decomposable\footnote{A non-decomposable element
  in $\mathfrak{B}_7^4$ can only involve at most $3$ distinct integers
  $\in \{1,2,3,4,5,6\}$. It can then be checked that, up to
  permutation and up to multiplication by an element in
  $(\Z/7)^\times$, the only non-decomposable elements in
  $\mathfrak{B}_7^4$ are $(1,1,2,2,4,4), (1,1,1,3,3,5)$ and
  $(1,1,1,1,5,5)$. But these are quasi-decomposable as one can see
  after adding the fourth and fifth digits. Also, one can check that
  the multiples of these three elements are either decomposable or
  quasi-decomposable. Finally, notice for example that $(2,2,2,2,2,4)$
  does not belong to $\mathfrak{B}_7^4$ because it has a multiple,
  namely $(1,1,1,1,1,2)$, whose digits add to $7$.}.  Hence, $
(H^{1,3}(X) \oplus H^{2,2}(X) \oplus H^{3,1}(X)) \cap H^4_{prim}(X,\Q)
\subseteq \N^1H^4(X)$.  Since $H^4(X)$ is the orthogonal sum (for the
cup product) of $H^4_{prim}(X,\Q)$ with the span of a linear section,
the result follows. \end{proof}

Let $G = \mu_7$ be the group of complex $7$-th roots of unity and let
$\zeta$ be a generator of $\mu_7$. We let $G$ act on $\P_\C^5$ in the
following way:
$$\zeta\cdot [x_0, x_1, x_2, x_3, x_4, x_5] =  [x_0, \zeta x_1, \zeta^2
x_2, \zeta^3 x_3, \zeta^4 x_4, \zeta^5 x_5].$$ Such an action
restricts to a free action of $G$ on $X$. Thus, the quotient $f : X \r
Y=X/G$ exists and defines a smooth projective $4$-fold $Y$. Since $X$
is of general type, so is $Y$.

\begin{lemma} The variety $Y$ above satisfies $H^i(Y,O_Y) = 0$ for $1
  \leq i \leq 4$.
\end{lemma}
\noindent \emph{Proof.} We have $H^i(Y,O_Y) = H^i(X,O_X)^G$. The weak
Lefschetz theorem settles the cases $1 \leq i \leq 3$. For the case
$i=4$, by Hodge theory we have $H^4(X,O_X) = \overline{H^0(X,K_X)}$.
Let $\Omega = \sum_i (-1)^i x_i dx_0 \wedge \ldots \widehat{dx_i}
\ldots \wedge dx_4$ be a generator of $H^0(\P^5,K_{\P^5}(6))$. We have
$$\zeta^*\Omega = \zeta\Omega.$$ A global
$4$-form $\omega \in H^0(X,K_X)$ on $X$ is the restriction of
$P\Omega/F$ to $X$ for some $P \in H^0(\P^5,O_{\P^5}(1))$ where $F =
x_0^7 + x_1^7 + x_2^7 +x_3^7 + x_4^7 +x_5^7\in H^0(\P^5,O_{\P^5}(7))$.
We have $$\zeta^*\omega = \zeta \Omega/F \cdot (\zeta^*P)|_X.$$
Now, $\zeta^*$ does not admit the eigenvalue $\zeta^{-1}$ on $
H^0(\P^5,O_{\P^5}(1))$. Hence, $H^i(X,O_X)^G=0$. \qed

\begin{proposition} The variety $Y$ above satisfies the Hodge
  conjecture as well as the Lefschetz standard conjecture and
  $\N^1H_4(Y) = H_4(Y)$. Moreover, $Y$ is dominated by a product of
  curves.
\end{proposition}
\noindent \emph{Proof.}  Let $p \in \CH_4(X \times X)$ denote the
correspondence $\frac{1}{7} {}^t\Gamma_{f} \circ \Gamma_{f}$. Because
$\Gamma_{f} \circ {}^t\Gamma_{f} = 7\Delta_Y \in \CH_4(Y \times Y)$,
$p$ defines an idempotent that identifies $H^4(Y)$ with $p_*H^4(X)$.
It follows that $Y$ satisfies the Hodge conjecture and that
$\N^1H_4(Y) = H_4(Y)$. The two other properties are straightforward.
\qed

\begin{theorem} \label{example} The variety $Y$ satisfies the
  following: $\CH^\Z_0(Y) = \CH^\Z_3(Y)= \Z$, $\CH_2(Y)$ is a finite
  dimensional $\Q$-vector space. Moreover, there exist a surface $S$
  and a correspondence $\Gamma \in \CH_3(S \times Y)$ such that
  $\CH_1(Y) = \Gamma_*\CH_0(S)$.
\end{theorem}
\noindent \emph{Proof.} The variety $Y$ satisfies $(\star)$ and
$(\star \star)$. We can therefore apply the results of the previous
sections : Proposition \ref{integralsupport} for $\CH^\Z_0(Y)$,
Proposition \ref{gen} for $\CH_1(Y)$ and Proposition \ref{aaa} for
$\CH_2(Y)$. The statement for $ \CH^\Z_3(Y)$ is obvious. \qed

\begin{remark}
  More precisely the Chow motive of $Y$ can be shown to ``come from''
  the Chow motive of a surface. Also, Murre's conjectures hold for $Y$
  (see Corollary \ref{FermatMurre}).
\end{remark}

\begin{remark} It can be shown that $H^3(Y,\Omega_Y) \neq 0$ and hence
  that $\Gr_{\N}^1H_4(Y) \neq 0$. By Lewis' and Schoen's theorem
  \ref{LS}, there is no curve $C$ and no correspondence $\Gamma \in
  \CH_2(C \times Y)$ such that $\CH_1(Y) = \Gamma_*\CH_0(C)$.
\end{remark}

\begin{remark}
  Unfortunately Ran's and Shioda's technique for proving the
  generalized Hodge conjecture for Fermat varieties does not work for
  high-degree Fermat varieties. For example, it is mentioned in
  \cite[\S 3 (13)]{Shioda2} that the method fails for Fermat $3$-folds
  of degree $11$ because there is an element in $\mathfrak{B}_{11}^3 =
  \{(b_0,b_1,b_2,b_3,b_4) \ | \ 1 \leq b_i \leq 10, \ \sum_{i=0}^4
  <tb_i>= 22 \ \mathrm{or} \ 33 \ \mathrm{for} \ \mathrm{all} \ t \in
  (\Z/11)^\times \}$ which is neither decomposable nor
  quasi-decomposable, namely $(1,1,5,7,8)$.
\end{remark}

\section{Niveau and Coniveau filtrations on Chow groups}

We use the construction of the idempotents $\Pi_{i,j}$ of Theorem
\ref{action} to relate some conjectures on algebraic cycles. Andr\'e
and Kahn \cite{AnKa} proved that the conjunction of Grothendieck's
standard conjectures with the Bloch--Beilinson--Murre conjectures imply
Kimura's finite-dimensionality conjecture. Jannsen \cite[Cor.
6.4(c)]{Jannsen4} proved that the conjunction of such conjectures
implies that the niveau and the coniveau filtrations on Chow groups
coincide. Here, assuming Grothendieck's standard conjectures and
Kimura's conjecture, we introduce a filtration on Chow groups using
our idempotents and, in the spirit of Jannsen's paper, prove that if
our filtration on Chow groups coincides with the coniveau filtration
considered by Jannsen then the BBM conjectures hold.

\paragraph{Murre's conjectures.} These were formulated in \cite{Murre1}
and were shown to be equivalent to those formulated by Bloch and
Beilinson in \cite[Th 5.2]{Jannsen}. Let $X$ be a smooth projective
variety of dimension $d$ over $k$. Murre conjectured that: \medskip

(A) $X$ has a Chow--K\"unneth decompositio: There exist mutually
orthogonal idempotents $\Pi_0, \ldots, \Pi_{2d} \in \CH_d(X \times X)$
adding to the identity such that $(\Pi_i)_*H_*(X)=H_i(X)$ for all $i$.

(B)  $\Pi_0, \ldots, \Pi_{2l-1},\Pi_{d+l+1}, \ldots, \Pi_{2d}$ act
trivially on $\CH_l(X)$ for all $l$.

(B') $\Pi_0, \Pi_1, \ldots, \Pi_{2l-1}$ act trivially on $\CH_l(X)$ for
all $l$.

(C) $F^i\CH_l(X) := \ker(\Pi_{2l}) \cap \ldots \cap \ker(\Pi_{2l+i-1})$
doesn't depend on the choice of the idempotents $\Pi_j$. Here the
idempotents $\Pi_j$ are acting on $\CH_l(X)$.

(D) $F^1\CH_l(X) = \CH_l(X)_\hom$. \medskip

Let $X$ be a smooth projective variety that satisfies $(\star \star)$
and whose K\"unneth projectors $\pi_i : H_*(X) \r H_i(X) \r H_*(X)$
are algebraic. By Lemma \ref{liftproj}, the projectors $\pi_i$ lift to
idempotents $\Pi_i \in \CH_d(X\times X)$ and hence $X$ satisfies
Murre's conjecture (A). Moreover if $(\Pi'_i)_{0 \leq i \leq 2d}$ is
another choice of liftings for the K\"unneth projectors then, still by
Lemma \ref{liftproj}, there exists a nilpotent correspondence $N \in
\CH_d(X \times X)$ such that $\ker(\Pi'_{2l}) \cap \ldots \cap
\ker(\Pi'_{2l+i-1}) = (1+N)_* \big( \ker(\Pi_{2l}) \cap \ldots \cap
\ker(\Pi_{2l+i-1})\big)$. We wish to show that if $X$ and $X \times X$
satisfy Murre's conjecture (B') then Murre's conjecture (C) holds for
$X$. For this purpose it is enough to show that $F^i\CH_l(X)$ is
stabilized by the action of correspondences in $\CH_d(X \times X)$.
The transpose idempotents $\widetilde{\Pi}_i := {}^t\Pi_{2d-i}$ define
another set of mutually orthogonal idempotents lifting the K\"unneth
projectors. The correspondences $\widetilde{\Pi}_i \times \Pi_j \in
\CH_{2d}((X\times X) \times (X \times X))$ then define mutually
orthogonal idempotents and $X\times X$ satisfies Murre's conjecture
(A) with respect to the idempotents $\Pi_k^{X\times X} :=
\sum_{i+j=k}\widetilde{\Pi}_i \times \Pi_j$ by the K\"unneth formula
in cohomology. The proof of the following proposition is essentially
contained in the proof of Jannsen's \cite[Th 5.2 \& Prop.
5.8]{Jannsen}.

\begin{proposition} \label{Murre-C}
  If $X$ satisfies $(\star \star)$ as well as Murre's conjectures (A)
  and (B'), and if $X \times X$ given with the idempotents
  $\Pi_k^{X\times X}$ above satisfies Murre's conjecture (B') then $X$
  satisfies Murre's conjecture (C).
\end{proposition}
\begin{proof} By the above discussion it is enough to prove that
  $F^i\CH_l(X)$ is stabilized by the action of correspondences $\alpha
  \in \CH_d(X \times X)$.  Because $X$ satisfies Murre's conjecture
  (B'), we have $F^i\CH_l(X) = \im(\Pi_{2l+i} + \ldots + \Pi_{2d})$. By
  orthogonality of the idempotents $\Pi_k$, it is thus equivalent to
  prove that $\Pi_j \circ \alpha \circ \Pi_k = 0$ for $2l+i \leq k
  \leq 2d$ and $2l \leq j \leq 2l+i-1$. But $\Pi_j \circ \alpha \circ
  \Pi_k = ({}^t\Pi_k \times \Pi_j)_*\alpha = (\widetilde{\Pi}_{2d-k}
  \times \Pi_j)_*\alpha$ and hence $\Pi_j \circ \alpha \circ \Pi_k \in
  \im(\Pi_{2d-k+j}^{X\times X})$. Now $2d-k+j < 2d$ and because $X
  \times X$ is assumed to satisfy Murre's conjecture (B') we get that
  $\Pi_{2d-k+j}^{X\times X}$ acts as zero on $\CH_d(X \times X)$.
\end{proof}

\begin{corollary} \label{2} If $(\star \star)$ and Murre's conjectures
  (A) and (B') hold for all smooth projective varieties then Murre's
  conjecture (C) holds also for all smooth projective varieties.
\end{corollary}

\begin{remark} \label{C'} Some authors replace Murre's conjecture (C)
  by the following weaker conjecture which is satisfied when $X$ has a
  CK decomposition and satisfies $(\star\star)$. \medskip

  (C') If $\{\Pi_i \}$ and $\{\Pi_i' \}$ are two Chow--K\"unneth
  decompositions for $X$ inducing, respectively, filtrations $F$ and
  $F'$ on the Chow groups of $X$, then there exists a nilpotent
  correspondence $N \in \CH_d(X \times X)$ such that
  $F^i\CH_l(X)=(1+N)_*(F')^i\CH_l(X)$ for all $i$ and all $l$. \medskip

  \noindent If $S$ is a smooth projective surface, then Murre
  \cite{Murre} proved that $S$ satisfies Murre's conjectures (A), (B)
  and (D). If, moreover, $S$ satisfies $(\star\star)$, then $S$ does
  actually satisfy Murre's conjecture (C) (and not just (C')).  This
  can be proved using Proposition \ref{Murre-C} or can be seen to be a
  direct consequence of Theorem \ref{Beauville2} to come.
\end{remark}

If $X$ satisfies property B, then the correspondences $\pi_i$ exist by
Kleiman \cite{Kleiman} and factor through a smooth projective variety
of dimension $\min(i,2d-i) $. Therefore, thanks to Lemma
\ref{supportlemma}, if $X$ satisfies property B and $(\star \star)$,
then the idempotents $\Pi_{d-l+1}, \ldots, \Pi_{2d}$ act trivially on
$\CH_l(X)$ for all $l$. Thus Murre's conjecture (B') is equivalent to
Murre's conjecture (B) for those varieties that satisfy B and $(\star
\star)$. \medskip

If $X$ satisfies $(\star)$ and $(\star \star)$, then we let $\Pi_i :=
\sum_j \Pi_{i,j}$. Theorem \ref{action} immediately gives the
following proposition.

\begin{proposition} \label{murre} If $X$ is a smooth projective
  variety that satisfies $(\star)$ and $(\star \star)$ then
    \begin{itemize}
  \item  $X$
  satisfies Murre's conjectures (B) if and only if
  $(\Pi_{i,j})_*\CH_l(X) = 0$ whenever $i < 2l$ and $i- j > l + 1$.
  \item  $X$
  satisfies Murre's conjectures (D) if and only if
  $(\Pi_{i,j})_*\CH_l(X) = 0$ whenever $i = 2l$ and $i- j > l + 1$.
  \end{itemize}
\end{proposition}

\paragraph{Niveau filtration on Chow groups.} Nori \cite{Nori}
introduced an increasing filtration on $\CH_l(X)_\hom$ :
$$\CH_l(X)_\alg = A_0\CH_l(X)_\hom \subseteq  A_1\CH_l(X)_\hom \subseteq
\ldots \subseteq A_l\CH_l(X)_\hom = \CH_l(X)_\hom$$ called the niveau
filtration and defined as follows :
$$A_r\CH_l(X)_\hom := \sum \im(\Gamma_* : \CH_r(Y)_\hom \r \CH_l(X)),$$
where the sum runs through all varieties $Y$ and all correspondences
$\Gamma \in \CH^{d-l+r}(Y \times X)$. From the definition it is clear
that this filtration is functorial with respect to the action of
correspondences, i.e. if $\alpha \in \CH_{d+k}(X \times Z)$ then
$\alpha_*A_r\CH_l(X)_\hom \subseteq A_r\CH_{l+k}(Z)_\hom$.\medskip

If $X$ satisfies $(\star)$ and $(\star \star)$, then we can define a
filtration $$\N_r\CH_l(X)_\hom := \big(\sum_i \sum_{j \geq l-r}
\Pi_{i,j}\big)_*\CH_l(X)_\hom.$$ By Theorem \ref{Ksupport} and Lemma
\ref{supportlemma} we have $\N_r\CH_l(X)_\hom \subseteq A_r\CH_l(X)_\hom
$.  If one assumes property B as well as $(\star \star)$ being valid
for all smooth projective varieties then this filtration is defined
for all smooth projective varieties.  Unfortunately, I cannot prove
(even if I assume Murre's conjectures) that this filtration is
functorial with respect to the action of correspondences. If this were
the case then the filtration would not depend on the choice of the
idempotents $\Pi_{i,j}$ as Lemma \ref{liftproj} shows.

\begin{proposition}
  Assume that property B and $(\star \star)$ hold for all smooth projective
  varieties. Then the filtration $\N$ is functorial if and only if it
  coincides with Nori's filtration $A$.
\end{proposition}

\begin{proof}
  The ``if'' part is obvious. Let's prove the ``only if'' part. Let
  $Y$ be a smooth projective variety and $\Gamma \in
  \CH^{d-l+r}(Y\times X)$. We have $\CH_r(Y)_\hom = \N_r\CH_r(Y)_\hom$.
  Therefore if the filtration $\N$ is functorial, then
  $\Gamma_*\CH_r(Y)_\hom \subseteq \N_r\CH_l(X)_\hom$. This proves $
  A_r\CH_l(X)_\hom \subseteq \N_r\CH_l(X)_\hom$.
\end{proof}

\paragraph{Coniveau filtration on Chow groups.} In
\cite[5.10(b)]{Jannsen4}, Jannsen introduces a coniveau filtration on
Chow groups (beware that in \textit{loc. cit.} this filtration is
denoted by $\N^\bullet)$
$$N^r\CH_l(X)_\hom = \sum \im(\Gamma_* : \CH^r(Y)_\hom \r
\CH_l(X)),$$ where the sum runs through all smooth projective varieties
$Y$ and all correspondences $\Gamma \in \CH_{r+l}(Y \times X)$. By
\cite[Prop. 5.3]{Nori}, we have $A_{r}\CH_l(X) \subseteq
N^{r+1}\CH_l(X)$. Jannsen proves the following theorem.

\begin{theorem}[Cor 6.4(c) \cite{Jannsen4}]
  Assume that Murre's conjectures as well as property B hold for all
  smooth projective varieties. Then $A_{r}\CH_l(X)_\hom =
  N^{r+1}\CH_l(X)_\hom$ for all smooth projective varieties $X$ and all
  integers $r$ and $l$.
\end{theorem}

Andr\'e and Kahn \cite{AnKa} showed that Murre's conjectures together
with Grothendieck's standard conjecture B imply Kimura's finite
dimensionality conjecture and thus that property $(\star \star)$ holds
for all smooth projective varieties. Thus, we prove a weak converse
(``weak'' because we cannot prove that the filtrations $\N$ and $A$
agree) to Jannsen's theorem.

\begin{proposition} \label{1} Let $X$ be a smooth projective variety
  that satisfies $(\star)$ and $(\star \star)$. If $\N_{r}\CH_l(X)_\hom
  = N^{r+1}\CH_l(X)_\hom$ for all all integers $r$ and $l$, then
  Murre's conjectures (A), (B) and (D) hold for $X$.
\end{proposition}
\begin{proof} Recall that by part (5) of Theorem \ref{action} we have $
  (\Pi_{i,j})_*\CH_l(X) \subseteq \CH_l(X)_\hom$ unless $i=2l$ and
  $j=l$.  Given $l$ and $r$, by Theorem \ref{Ksupport} and Lemma
  \ref{supportlemma} we have $$\big(\sum_{l < i-j \leq l+r} \Pi_{i,j}
  \big)_*\CH_l(X) \subseteq N^r\CH_l(X)_\hom.$$ Let's now be given $l$.
  By Theorem \ref{action} we already know that $\Pi_{i,j}$ acts
  trivially on $\CH_l(X)$ when $i\leq 2l$ and $i-j=l+1$. Let's suppose
  that we have proved that $\Pi_{i,j}$ acts trivially on $\CH_l(X)$
  when $i\leq 2l$ and $l < i-j \leq l+r$ for some $r \geq 1$. Because
  we are assuming $\N_{r}\CH_l(X)_\hom = N^{r+1}\CH_l(X)_\hom$ we get
  that $\Pi_{i,j}$ acts trivially on $\CH_l(X)$ when $i-j = l+r+1$ and
  $j<l-r$, that is when $i-j = l+r+1$ and $i<2l+1$. Therefore, by
  induction on $r$, this holds for all positive integers $r$.  This
  proves that $\Pi_{i,j}$ acts trivially on $\CH_l(X)$ whenever $i-j >
  l$ and $i \leq 2l$.  Proposition \ref{murre} then shows that $X$
  satisfies Murre's conjectures (B) and (D).
\end{proof}

Combining Proposition \ref{1} and Corollary \ref{2} gives the
following proposition.

\begin{proposition} \label{end} Assume that properties B and $(\star
  \star)$ hold for all smooth projective varieties. If
  $\N_{r}\CH_l(X)_\hom = N^{r+1}\CH_l(X)_\hom$ for all smooth projective
  varieties $X$ and all integers $r$ and $l$, then Murre's conjectures
  hold for all smooth projective varieties. \qed
\end{proposition}

\section{Murre's conjectures in some new cases}
\label{exmurre}

A prerequisite for proving Murre's conjectures for a given variety $X$
is to have a Chow--K\"unneth decomposition for $X$ at our disposal.
Either we know that $X$ satisfies properties B and $(\star \star)$, in
which case the work done in \S1 may apply; or $X$ is given with a
specific CK decomposition. In the first paragraph, we go through the
known cases of varieties satisfying property B and also through the
known cases of varieties having a CK decomposition. In the second
paragraph, we use Theorem \ref{action} to establish Murre's
conjectures for some varieties satisfying $(\star\star)$. In the third
paragraph, we establish Murre's conjectures for some varieties having
a CK decomposition but for which we cannot prove that they satisfy
$(\star\star)$. Finally, in the fourth paragraph, we illustrate our
results with explicit examples.

\subsection{Preliminaries}

\paragraph{Varieties satisfying property B.} Let $E_B$ be the set of
smooth projective varieties satisfying property B. It is known that
curves, surfaces, complete intersections and abelian varieties belong
to $E_B$. We have the following lemma.

\begin{lemma} \label{Bdominant}
  The set $E_B$ is stable under product and smooth hypersurface
  section. If $X \in E_B$ and if $f : X \r Y$ is a dominant morphism,
  then $Y \in E_B$. If $X \in E_B$ and if $Z$ is a smooth subvariety
  of $X$ that belongs to $E_B$, then the blow-up of $X$ along $Z$
  belongs to $E_B$.
\end{lemma}

\begin{proof}
  Stability under product and hypersurface section is well known. The
  rest is contained in \cite[Lemma 4.2]{Arapura}.
\end{proof}

\paragraph{Varieties having a Kimura finite-dimensional Chow motive.}
Let $E_K$ be the set of smooth projective varieties having a Kimura
finite-dimensional Chow motive \cite{Kimura}. Curves, abelian
varieties and Fermat hypersurfaces \cite{KS} belong to $E_K$.

\begin{lemma} \label{Kdominant} The set $E_K$ is stable under product.
  If $X \in E_K$ and if $f : X \r Y$ is a dominant morphism, then $Y
  \in E_K$. If $X \in E_K$ and if $Z$ is a smooth subvariety of $X$
  that belongs to $E_K$, then the blow-up of $X$ along $Z$ belongs to
  $E_K$. Moreover if $X \in E_K$ then $X$ satisfies $(\star\star)$.
\end{lemma}
\begin{proof}
  For the proof, see \cite{Kimura}.
\end{proof}

 \paragraph{Varieties admitting a CK decomposition.} Let $E_{CK}$ be
 the set of smooth projective varieties having a CK decomposition. As
 for $E_B$, it is known that curves, surfaces, complete intersections
 and abelian varieties belong to $E_{CK}$.

 \begin{lemma} \label{CKdominant} The set $E_{CK}$ is stable under
   product.   If $X \in E_{CK}$ and if $Z$ is a
   smooth subvariety of $X$ that belongs to $E_{CK}$, then the
   blow-up of $X$ along $Z$ belongs to $E_{CK}$.
\end{lemma}
\begin{proof}
  The first point is straightforward and the second point follows from
  the blow-up formula for Chow motives.
\end{proof}

\begin{remark}
  It is tempting to think that if $Y$ is a variety dominated by $X \in
  E_{CK}$, then $Y \in E_{CK}$. Indeed, if $\{\Pi_i\}$ is a CK
  decomposition for $X$ and if $f : X \r Y$ is a dominant (generically
  finite for convenience) map, then $\{\frac{1}{\deg f} \Gamma_f \circ
  \Pi_i \circ {}^t\Gamma_f\}$ gives a K\"unneth decomposition for $X$
  (this is because the action of $\Pi_i$ on homology is central) and
  it is tempting to think that they should give a CK decomposition for
  $Y$. This last point is, however, far from being clear.
\end{remark}

\subsection{Murre's conjectures using Theorem \ref{action}}

\paragraph{4.2.1.} We start with an immediate consequence of Theorem
\ref{action}. Recall that a fourfold that satisfies property B
satisfies $(\star)$ too.

\begin{theorem} \label{Beauville}
  If $X$ is a fourfold that satisfies B and $(\star \star)$, then $X$
  satisfies Murre's conjectures (A), (B) and  (C').
\end{theorem}

Beauville \cite{Beauville} had already proved that result for abelian
fourfolds using a different technique. Because abelian fourfolds
satisfy both  B and $(\star \star)$, Theorem \ref{Beauville} can be
seen as a generalization of Beauville's result. \medskip

Note that the only obstruction for $X$ to satisfying Murre's
conjecture (D) is due to the fact that we cannot prove that $\Pi_{4,0}$
acts trivially on $\CH_2(X)$ if $H_4(X) \neq \N^1H_4(X)$.

\paragraph{4.2.2.} Let $F$ be the set of Kimura finite-dimensional
smooth projective varieties $X$ that satisfy B and whose cohomology
satisfies $H_i(X) = \N^{\lfloor i/2 \rfloor -1}H_i(X)$ for all $i$. In
other words, the cohomology of $X \in F$ in even degree is generated,
via the action of correspondences, by the degree $2$ homology of
surfaces and the cohomology of $X$ in odd degrees is generated by the
degree $3$ homology of $3$-folds. In particular, $X \in F$ satisfies
property $(\star)$ and $(\star \star)$. Let's first mention that the
set $F$ is not too small.

\begin{proposition} \label{stabF} The set $F$ contains Kimura finite
  dimensional $3$-folds that satisfy B (e.g.  $3$-folds dominated by a
  product of curves, abelian $3$-folds, Fermat $3$-folds and
  rationally connected $3$-folds).  If $X \in F$ and if $Z$ is a
  smooth subvariety of $X$ of dimension $\leq 3$ that satisfies B and
  $(\star \star)$, then the blow-up of $X$ along $Z$ belongs to $F$.
  If $X \in F$ and if $f : X \r Y$ is a dominant morphism, then $Y \in
  F$.
\end{proposition}

\begin{proof}
  This follows from Lemmas \ref{Bdominant} and \ref{Kdominant}.
\end{proof}

\begin{corollary} \label{Fdominant}
  Let $X$ be a $3$-fold belonging to $F$ and let $f : X
  \dashrightarrow Y$ be a dominant rational map. Then $Y$ belongs to $F$.
\end{corollary}

\begin{proof}
  By resolution of singularities, there exists a sequence of blow-ups
  $\widetilde{X}_n \r \ldots \r \widetilde{X}_1 \r X$ along smooth
  curves and a dominant map $\widetilde{X}_n \r Y$. By Proposition
  \ref{stabF}, $\widetilde{X}_n \in F$ and hence $Y \in F$.
\end{proof}

Once again, as an immediate consequence of Theorem \ref{action}, we
get the following theorem.

\begin{theorem} \label{Beauville2} If $X \in F$, then $X$ satisfies
  Murre's conjectures (A), (B), (C') and (D).

  \noindent If, moreover, $H_{2i+1}(X) = \N^{i}H_{2i+1}(X)$ for all $i$
  (for example if $X$ is a surface that satisfies $(\star \star)$),
  then $X$ satisfies Murre's conjecture (C).
\end{theorem}

\begin{proof}
  Let's only prove the second point. If $X \in F$ with $H_{2i+1}(X) =
  \N^{i}H_{2i+1}(X)$ for all $i$, then, for all $l$, Proposition
  \ref{gen} provides the existence of a smooth projective surface $S$,
  of a correspondence $\Gamma \in \CH_{l+2}(S \times X)$ and of a
  correspondence $\Gamma' \in \CH_{d-l}(X \times S)$ such that $\Gamma
  \circ \Gamma'$ acts as the identity on $\CH_l(X)$. (This is the case
  $p=2$ and $q=0$ of the proposition). Because algebraic and
  homological equivalence agree on surfaces, we get that they also
  agree on $X$. Therefore, by Theorem \ref{action}, we have that
  $F^1\CH_l(X)= \CH_l(X)_\hom$, $F^2\CH_l(X)=\ker (AJ_l : \CH_l(X)_\hom \r
  J_l(X))$ and $F^3\CH_l(X)=0$ for all $l$. In particular, the
  filtration does not depend on the choice of a CK decomposition for
  $X$.
\end{proof}

\begin{corollary}
  If $X$ is rationally dominated by a product of three curves, then
  $X$ satisfies Murre's conjectures (A), (B), (C') and (D).
\end{corollary}

\begin{corollary} \label{C:abelian}
  Abelian $3$-folds satisfy Murre's conjectures (A), (B), (C') and
  (D).
\end{corollary}

\begin{remark} Murre's conjectures (B) and (D) were already known to
  hold for abelian $3$-folds by \cite{Beauville}.
\end{remark}

\begin{corollary}\label{FermatMurre}
  The variety $Y$ of \S \ref{Fermatexample} satisfies Murre's
  conjectures (A), (B), (C) and (D).
\end{corollary}
\begin{proof}
  The variety $Y$ is dominated by a product of curves and satisfies
  $H_3(Y) = H_5(Y) = 0$ and $H_4(Y) = \N^1H_4(Y)$. Therefore $Y$
  satisfies the assumptions of Theorem \ref{Beauville2}.
\end{proof}

\subsection{Murre's conjectures for some varieties having a CK
  decomposition}

Let $G$ be the set of smooth projective varieties $X$ that have a CK
decomposition $\{\Pi_i\}_{0 \leq i \leq 2\dim X}$ such that for all
$i$ there exist subvarieties $Y_i$ and $Z_i$ of respective dimension
$i+1$ and $i+2$ such that $\Pi_{2i}$ has a representative supported on
$X \times Y_i$ and $\Pi_{2i+1}$ has a representative supported on $X
\times Z_i$. The correspondence $\Pi_{2i}$ should be thought of as
factoring through a surface and the correspondence $\Pi_{2i+1}$ should
be thought of as factoring through a $3$-fold (although beware that
this type of condition is more restrictive).  Obvious examples of
varieties belonging to $G$ are given by curves, surfaces, products of
a curve and a surface, and smooth complete intersections of dimension
$3$. We will give more examples of varieties belonging to $G$ in the
next paragraph.  For the moment, let's only mention that the set $G$
is not too small:

\begin{proposition} \label{belongG} If $X \in G$ and if $Z$ is a
  smooth subvariety of $X$ that belongs to $G$, then the blow-up of
  $X$ along $Z$ belongs to $G$.
\end{proposition}
\begin{proof}
  This follows from the blow-up formula for Chow motives.
\end{proof}

\begin{remark}
  If $X$ is a variety that belongs to $F$, then its CK projectors
  $\Pi_{2i}$ (resp. $\Pi_{2i+1})$ are homologically equivalent to some
  cycles supported on $X\times Y_i$ (resp. $X\times Z_i)$ for some
  subvariety $Y_i$ of dimension $i+1$ (resp. $Z_i$ of dimension
  $i+2)$. However, although it is expected to be the case, it is not
  clear that the projectors themselves are supported on some $X\times
  Y_i$ or $X\times Z_i$. Thus it is not clear that $F \subseteq G$.
  The reverse inclusion is even less clear but $F=G$ would follow from
  general conjectures on algebraic cycles.
\end{remark}

The main result of this section is the following theorem.

\begin{theorem} \label{E}
  If $X \in G$, then $X$ satisfies Murre's conjectures (A), (B') and
  (D). If, moreover, $X$ satisfies $(\star \star)$, then $X$ satisfies
  Murre's conjecture (C').
\end{theorem}

We split its proof into two lemmas. (the statement about conjecture
(C') follows from Remark \ref{C'}).

The first lemma settles Murre's conjecture (D) for varieties
belonging to $G$.

\begin{lemma} \label{D}
  If $X \in G$, then $\ker \big( \Pi_{2i} : \CH_i(X) \r \CH_i(X)\big) =
  \CH_i(X)_\hom$ for all $i$.
\end{lemma}
\begin{proof}
  Let's fix $i$. The functoriality of the cycle class map to singular
  homology with respect to the action of correspondences combined with
  the fact that $\Pi_{2i}$ acts as the identity on $H_{2i}(X)$
  immediately implies that $\ker (\Pi_{2i}) \subseteq \CH_i(X)_\hom$.

  Let's now consider $\widetilde{Y}_i$ a desingularization of $Y_i$.
  Then the idempotent $\Pi_{2i}$ factors through $\widetilde{Y}_i$,
  i.e.  there exist $g \in \CH_{d}(X \times \widetilde{Y}_i)$ and $f
  \in \CH_{i+1}(\widetilde{Y}_i \times X)$ such that $\Pi_{2i} = f
  \circ g$. In particular, the action of $\Pi_{2i}$ on $\CH_i(X)_\hom$
  factors through $\CH^1(\widetilde{Y}_i)_\hom$. By functoriality of
  the Abel--Jacobi map, we have the following commutative diagram :
  \begin{center} $ \xymatrix{ \CH_i(X)_\hom \ar[d]^{AJ_i} \ar[r]^{g_*}
      & \CH^1(\widetilde{Y}_i)_\hom \ar[d]^{\simeq} \ar[r]^{f_*} &
      \CH_i(X)_\hom
      \ar[d]^{AJ_i}  \\
      J_{i}(X) \ar[r] & \mathrm{Pic}^0(\widetilde{Y}_i) \ar[r] &
      J_{i}(X).}$ \end{center} The composite of the two bottom
  arrows is zero because $\Pi_{2i}$ acts trivially on $H_{2i+1}(X)$.
  Therefore, if $\alpha \in \CH_i(X)_\hom$, then we have $(g\circ
  f\circ g)_*\alpha = 0$ and hence $f_* \circ (g\circ f\circ
  g)_*\alpha = 0$, i.e. $(\Pi_{2i} \circ \Pi_{2i})_*\alpha = 0$, that
  is, $(\Pi_{2i})_*\alpha = 0$.
\end{proof}

The second lemma settles Murre's conjecture (B') for varieties
belonging to $G$.

\begin{lemma} \label{B'}
  If $X \in G$, then $\Pi_i$ acts trivially on $\CH_l(X)$ for all $i <
  2l$.
\end{lemma}
\begin{proof}
  For obvious dimension reasons it is enough to prove that $\Pi_{2i}$
  acts trivially on $\CH_{i+1}(X)$ and that $\Pi_{2i+1}$ acts trivially
  on $\CH_{i+1}(X)$ and on $\CH_{i+2}(X)$.

  Let's consider $\widetilde{Z}_i$ a desingularization of $Z_i$.
  Then the idempotent $\Pi_{2i+1}$ factors through $\widetilde{Z}_i$,
  i.e.  there exist $g \in \CH_{d}(X \times \widetilde{Z}_i)$ and
  $f \in \CH_{i+2}(\widetilde{Z}_i \times X)$ such that $\Pi_{2i} = f
  \circ g$.

  Let's first prove that $\Pi_{2i+1}$ acts trivially on $\CH_{i+2}(X)$.
  By functoriality of the cycle class map, we have the commutative
  diagram \begin{center} $ \xymatrix{ \CH_{i+2}(X) \ar[d]^{cl_{i+2}}
      \ar[r]^{g_*} & \CH^0(\widetilde{Z}_i) \ar[d]^{\simeq}
      \ar[r]^{f_*} & \CH_{i+2}(X)  \ar[d]^{cl_{i+2}}  \\
      H_{2i+4}(X) \ar[r] & H^0(\widetilde{Z}_i) \ar[r] &
      H_{2i+4}(X).}$
  \end{center} By definition of a CK decomposition, the idempotent
  $\Pi_{2i+1}$ acts trivially on $H_j(X)$ for $j \neq 2i+1$. Therefore
  if $\alpha \in \CH_{i+2}(X)$, then $(g\circ f\circ g)_*\alpha = 0$
  and hence $f_* \circ (g\circ f\circ g)_*\alpha = 0$, i.e.
  $(\Pi_{2i+1} \circ \Pi_{2i+1})_*\alpha = 0$, that is
  $(\Pi_{2i+1})_*\alpha = 0$.

  The fact that $\Pi_{2i}$ acts trivially on $\CH_{i+1}(X)$ is similar.

  Let's now prove that $\Pi_{2i+1}$ acts trivially on $\CH_{i+1}(X)$.
  For this purpose, let's consider the cycle class map to Deligne
  cohomology. This map is functorial with respect to the action of
  correspondences and induces an isomorphism $cl_{\mathcal{D}}^1 :
  \CH^1(Y) \stackrel{\simeq}{\longrightarrow}
  H_{\mathcal{D}}^2(Y,\Q(1))$ for any smooth projective variety $Y$.
  Once again we have a commutative diagram \begin{center} $ \xymatrix{
      \CH_{i+1}(X) \ar[d]^{cl_{\mathcal{D}}} \ar[r]^{g_*} &
      \CH^1(\widetilde{Z}_i) \ar[d]^{\simeq}
      \ar[r]^{f_*} & \CH_{i+1}(X)  \ar[d]^{cl_{\mathcal{D}}}  \\
      H_{\mathcal{D}}^{2d-2i-2}(X,\Q(d-i-1)) \ar[r] &
      H^2_{\mathcal{D}}(\widetilde{Z}_i,\Q(1)) \ar[r] &
      H_{\mathcal{D}}^{2d-2i-2}(X,\Q(d-i-1)).}$
  \end{center}
  The idempotent $\Pi_{2i+1}$ acts trivially on $H_{2i+2}(X)$ and it
  acts  trivially also on the intermediate Jacobian $J_{i+1}(X)$
  because it acts trivially on $H_{2i+3}(X)$.  Therefore, $\Pi_{2i+1}$
  acts trivially on $H_{\mathcal{D}}^{2d-2i-2}(X,\Q(d-i-1))$ because
  this last group is an extension of the Hodge classes in
  $H_{2i+2}(X)$ by $J_{i+1}(X)$. Therefore, as before, if $\alpha \in
  \CH_{i+1}(X)$, then $(g\circ f\circ g)_*\alpha = 0$ and hence $f_*
  \circ (g\circ f\circ g)_*\alpha = 0$, i.e. $(\Pi_{2i+1} \circ
  \Pi_{2i+1})_*\alpha = 0$, that is, $(\Pi_{2i+1})_*\alpha = 0$.
\end{proof}

\begin{remark} Let $X$ be a smooth projective variety that has a CK
  decomposition $\{\Pi_i\}$. If $\Pi_{2i}$ has a representative
  supported on $X \times Y_i$ with $\dim Y_i = i+2$, then, using the
  same technique as in the two previous lemmas, it can be shown that
  $\Pi_{2i}$ acts trivially on $\CH_l(X)$ for $l > i$.  Likewise, it
  can be shown that if $\dim Y_i = i+n$ for some positive integer $n$,
  then $\Pi_{2i}$ acts trivially on $\CH_l(X)$ for $l > i+n-2$.
\end{remark}

\subsection{Examples}

\paragraph{4.4.1. Varieties belonging to $F$.} In \cite{Vial3}, it is
proved that if $X$ is a smooth projective variety with representable
Chow groups, then $X$ satisfies B, $H_i(X) = \N^{\lfloor i/2
  \rfloor}H_i(X)$ for all $i$ and the Chow motive of $X$ is finite
dimensional in the sense of Kimura. Examples of varieties with
representable Chow groups are discussed in \cite{Vial3} and include
curves, surfaces not of general type with vanishing geometric genus,
Godeaux surfaces, Barlow surfaces, rationally connected $3$-folds and
hypersurfaces of very low degree.

We immediately see that if $Y'$ is the product of three varieties with
representable Chow groups, then $Y'$ belongs to $F$. By Proposition
\ref{stabF}, a variety $Y$ obtained by repeatedly blowing up $Y'$
along smooth curves belongs to $F$. Still by Proposition \ref{stabF},
any variety dominated by $Y$ belongs to $F$. Therefore Theorem
\ref{Beauville2} gives the following theorem.

\begin{theorem} \label{rep} Let $X_1$, $X_2$ and $X_3$ be smooth
  projective varieties with representable Chow groups and let $Z$ be a
  variety dominated by a variety obtained by repeatedly blowing up
  $X_1 \times X_2 \times X_3$ along smooth curves. Then $Z$ satisfies
  Murre's conjectures (A), (B), (C') and (D).

  \noindent If, moreover, $\CH_l(X_3)$ is a finite-dimensional
  $\Q$-vector space for all $l$, then $Z$ satisfies Murre's conjecture
  (C).
\end{theorem}

\begin{proof}
  Only the last point deserves treatment. If $\CH_l(X_3)$ is a finite
  dimensional $\Q$-vector space for all $l$, then \cite[Th. 5]{Vial3}
  $H_{2i}(X_3) = \N^{i}H_{2i}(X_3)$ and $H_{2i+1}(X_3) = 0$ for all
  $i$, and the Chow motive of $X$ is finite-dimensional in the sense
  of Kimura.  As such, $Z$ belongs to $F$ and satisfies
  $H_{2i+1}(Z)=\N^iH_{2i+1}(Z)$ for all $i$. We are thus reduced to
  the last point of Theorem \ref{Beauville2}.
\end{proof}

\paragraph{4.4.2. Varieties belonging to $G$.} We wish now to give a
criterion on Chow groups for a variety to belong to $G$.

Before we proceed, let's define the following invariance property
which is weaker than Murre's conjecture (C).  Let $X$ be a smooth
projective variety that admits a CK decomposition $\{\Pi_k\}_{0 \leq k
  \leq 2d}$. We say that the CK decomposition $\{\Pi_k\}_{0 \leq k
  \leq 2d}$ is \textit{special} if, for all $i \neq d$, $\Pi_i$
factors through a $0$-dimensional variety if $i$ is even and $\Pi_i$
factors through a curve if $i$ is odd. By Murre \cite{Murre}, every
surface has a special CK decomposition.

Then we say that $X$ satisfies Murre's conjecture (C'') if \medskip

(C'') For any two special CK decompositions of $X$, the induced
filtrations on the Chow groups of $X$ coincide.

\begin{lemma} \label{special-iso}
  If $\{\Pi_k\}_{0\leq k \leq 2d}$ and $\{\Pi'_k\}_{0\leq k \leq 2d}$
  are special CK decompositions for $X$, then for every $i$ the Chow
  motives $(X, \Pi_k)$ and $(X,\Pi'_k)$ are isomorphic. Here $(X,
  \Pi_k)$ denotes the image of the idempotent endomorphism $\Pi_k$ of
  the Chow motive $h(X)$.
\end{lemma}
\begin{proof}
  When $k\neq d$ the proof is clear: both $(X,\Pi_k)$ and $(X,\Pi'_k)$
  are Kimura finite-dimensional, so that the composite of the
  embedding $(X,\Pi_k) \r h(X)$ with the projection $h(X) \r
  (X,\Pi'_k)$ is an isomorphism because it is an isomorphism modulo
  homological equivalence.  Similarly, if $M$ and $M'$ are the
  respective direct sums of the $(X,\Pi_k)$ and $(X,\Pi_k')$ for
  $k\neq d$, then the composite of the embedding $s : M \r h(X)$ with
  the projection $r' : h(X) \r M'$ is an isomorphism $u$. Thus if $r :
  h(X) \r M$ is the projection, then $r$ and $u^{-1} \circ r'$ are
  both left inverse to $s$, so that
  $$(X,\Pi_d) \simeq \ker r \simeq \mathrm{Coker} \, s \simeq \ker
  (u^{-1} \circ r') \simeq \ker r' \simeq (X,\Pi_d'),$$ as wanted.
 \end{proof}

Let $\{\Pi_k\}_{0\leq k \leq 2d}$ be a CK decomposition for $X$.  Call
$\{\Pi_k\}_{0\leq k \leq 2d}$ \emph{very special} if it is special and
if $\Pi_d$ has a representative supported on $X \times Y_n$ with $Y_n$
of dimension $n + 1$ when $d = 2n$ is even, and on $X \times Z_n$ with
$Z_n$ of dimension $n + 2$ when $d = 2n + 1$ is odd. One of the
reasons for introducing this notion is the following straightforward
lemma.

\begin{lemma} \label{CKG}
  If $X$ has a very special CK decomposition then $X$ belongs to the
  set $G$.
\end{lemma}

Let's rephrase in terms of motives what it means for $\{\Pi_k\}_{0
  \leq k \leq 2d}$ to be very special. In the odd-dimensional case for
example, the condition on the support is equivalent to requiring that
$\Pi_d$ factor through $h(\widetilde{Z}_n)$ for a resolution of
singularities $\widetilde{Z}_n \r Z_n$. Now an idempotent endomorphism
$e$ in a pseudo-abelian category factors through an object $N$ if and
only if the image of $e$ is a direct summand of $N$.  Thus the
condition on $\Pi_d$ is equivalent to requiring that $(X,\Pi_d)$ be a
direct summand of $h(\widetilde{Z}_n)$.

Note that if $\{\Pi_k\}_{0\leq k \leq 2d}$ is special, then by the
above lemma $(X,\Pi_i)^\vee$ (the dual of $(X,\Pi_i))$ is isomorphic
to $(X,\Pi_{2d-i},-d)$. Therefore, taking duals gives the further
equivalent conditions that $(X,\Pi_d)$ be a direct summand of
$h(\widetilde{Z}_n)(n-1)$ or that $\Pi_d$ factor through
$h(\widetilde{Z}_n)(n-1)$. Using again the above lemma, we also obtain
the following.

\begin{lemma} \label{veryspecial}
  If X has a very special CK decomposition then every special CK
  decomposition of X is very special.
\end{lemma}

We now give a criterion on the Chow groups of $X$ for $X$ to have a
very special CK decomposition. We start with the even-dimensional
case. The following is taken from \cite[Th.  4.5]{Vial3}.

\begin{theorem}
  Let $X$ be a smooth projective variety of even dimension $d=2n$. If
  $\CH_0(X)_\alg, \CH_1(X)_\alg, \ldots, \CH_{n-2}(X)_\alg$ are
  representable, then $X$ has a very special CK decomposition. \qed
\end{theorem}

\begin{proposition}
  Suppose that $d = 2n$ is even and that $X$ has a very special CK
  decomposition.  Then $X$ satisfies Murre's conjectures (A), (B),
  (C'') and (D). Precisely, algebraic and homological equivalence
  agree on $i$-cycles on $X$ for all $i$, and $AJ_i$ is injective for
  $i \neq n-1$. If $F$ is the filtration defined by a special CK
  decomposition for $X$, then $F^1\CH_i(X) = \CH_i(X)_{\mathrm{hom}}$
  for every $i$ and $F^2\CH_i(X) = 0$ for $i \neq n-1$, while
  $F^2\CH_{n-1}(X) = \ker(AJ_{n-1} : \CH_{n-1}(X)_{\mathrm{hom}} \r
  J_{n-1}(X))$ and $F^3\CH_{n-1}(X) = 0$.
\end{proposition}
\begin{proof} Let $\{\Pi_k\}_{0\leq k \leq 2d}$ be a special CK
  decomposition for $X$, which by Lemma \ref{veryspecial} is very
  special.  By Lemma \ref{CKG}, $X$ belongs to $G$ and Theorem \ref{E}
  shows that $X$ satisfies Murre's conjectures (A), (B') and (D) with
  respect to the CK decomposition $\{\Pi_k\}_{0\leq k \leq 2d}$.
  Let's now show that $X$ satisfies Murre's conjecture (B). This is
  clear because when $i \neq d$, $\Pi_i$ acts non-trivially only on
  $\CH_{\lfloor i/2 \rfloor}(X)$ because $\Pi_i$ factors through a
  curve and acts only in one degree in homology (see also \cite[Prop.
  2.9]{Vial3}); and when $i=d$, ${}^t\Pi_d = \Pi_d$ acts trivially on
  $\CH_l(X)$ for $l < n-1$ for dimension reasons (the action of
  $({}^t\Pi_d)_*$ on $\CH_l(X)$ factors through $\CH_{l+1-n}(Y_n)$ or
  the motive $(X,\Pi_d)$ is a direct summand of $h(\tilde{Y}_n)(n-1)$
  for a resolution of singularities $\tilde{Y}_n \r Y_n)$.

  The above shows that that $F^1\CH_i(X) = \CH_i(X)_\hom$ for all $i$.
  It also shows that $\CH_i(X) = (\Pi_{2i} + \Pi_{2i+1})_*\CH_i(X)$ for
  all $i \neq n-1$. Consequently, for all $i \neq n-1$, homological
  and algebraic equivalence agree on $\CH_i(X)$ and $F^2\CH_i(X)=0$.
  For dimension reasons, the idempotents $\Pi_i$ act trivially on
  $\CH_{n-1}(X)$ for $i < d-2$ and for $i > d$.  Therefore,
  $F^3\CH_{n-1}(X)=0$.  Moreover, the action of $\Pi_i$ on
  $\CH_{n-1}(X)$ factors through the Chow group of zero-cycles of a
  zero-dimensional variety if $i=d-2$, of a curve if $i=d-1$ and of
  $Y_n$ if $i=d$.  Therefore, homological and algebraic equivalence
  agree also on $\CH_{n-1}(X)$. By Theorem \ref{Chowproj} (see also
  \cite[Prop. 2.10]{Vial3}), we have that $F^2\CH_{n-1}(X) =
  \ker(AJ_{n-1} : \CH_{n-1}(X)_{\mathrm{hom}} \r J_{n-1}(X))$.

  We have thus showed that the filtration $F$ on the Chow groups of
  $X$ induced by the special CK decomposition $\{\Pi_k\}_{0\leq k \leq
    2d}$ does not depend on the choice of $\{\Pi_k\}_{0\leq k \leq
    2d}$. This settles (C'') for $X$.
\end{proof}

Combining the results above, we immediately obtain the following theorem.

\begin{theorem} \label{pair} Let $X$ be a smooth projective variety of
  even dimension $d=2n$. If $\CH_0(X)_\alg, \CH_1(X)_\alg, \ldots,
  \CH_{n-2}(X)_\alg$ are representable, then $X$ satisfies Murre's
  conjectures (A), (B), (C'') and (D). \qed
\end{theorem}

\begin{remark}
  When $n=1$, $X$ is a surface and no assumption is made on the Chow
  groups of $X$ in the theorem above.  In this case, conjectures (A),
  (B) and (D) were settled by Murre \cite{Murre1}.
\end{remark}

\begin{corollary} Let $X$ be a smooth projective fourfold.  If either
  $X$ is rationally connected or if $X$ admits a curve $C$ as a base
  for its maximal rationally connected fibration, i.e if there exists
  a rational map $f : X \dashrightarrow C $ with rationally connected
  general fiber, then $X$ satisfies Murre's conjectures (A), (B),
  (C'') and (D).
\end{corollary}
\begin{proof}
  Indeed if $X$ is rationally connected, then $\CH_0(X)=\Q$. If $X$
  admits a curve $C$ as a base for its maximal rationally connected
  fibration, then $\CH_0(X)_\alg$ is representable.
\end{proof}

\begin{examples}
  All smooth projective varieties that are birational to a Fano
  fourfold are rationally connected by Koll\'ar \cite{Kollar}, and hence
  satisfy Murre's conjectures (A), (B), (C'') and (D).
 \end{examples}
 \begin{examples}
   Let $Y$ be a generic hypersurface of $\P^4 \times \P^n$ of bidegree
   $(a,b)$ with $a \leq 4$. The projection of $Y$ onto $\P^n$ has
   rationally connected general fiber of dimension $3$. Let $C$ be a
   generic curve in $\P^n$. If $X := Y \times_{\P_n} C$, then the
   projection $X\r C$ has rationally connected general fiber.
   Therefore, $X$ satisfies Murre's conjectures (A), (B), (C'') and
   (D).
\end{examples}

Let's now deal with the odd-dimensional case.  In
\cite[Th.4.10]{Vial3} we proved the following.

\begin{theorem}
  Let $X$ be a smooth projective variety of odd dimension $d=2n+1$
  with $H^{n+1}(X,\Omega_X^{n-1})=0$. If $\CH_0(X)_\alg, \CH_1(X)_\alg,
  \ldots, \CH_{n-2}(X)_\alg$ are representable, then $X$ has a very
  special CK decomposition. \qed
\end{theorem}

\begin{proposition}
  Suppose that $d = 2n+1$ is odd and that $X$ has a very special CK
  decomposition.  Then $X$ satisfies Murre's conjectures (A), (B),
  (C'') and (D). Precisely, algebraic and homological equivalence
  agree on $i$-cycles on $X$ for $i \neq n$, and $AJ_i$ is injective
  for $i \neq n-1$. If $F$ is the filtration defined by a special CK
  decomposition for $X$, then $F^1\CH_i(X) = \CH_i(X)_{\mathrm{hom}}$
  for every $i$ and $F^2\CH_i(X) = 0$ for $i \neq n-1$, while
  $F^2\CH_{n-1}(X) = F^3\CH_{n-1}(X) = \ker(AJ_{n-1} :
  \CH_{n-1}(X)_{\mathrm{hom}} \r J_{n-1}(X))$ and $F^4\CH_{n-1}(X) = 0$.
\end{proposition}
\begin{proof} Let $\{\Pi_k\}_{0\leq k \leq 2d}$ be a special CK
  decomposition for $X$, which by Lemma \ref{veryspecial} is very
  special.  By Lemma \ref{CKG}, $X$ belongs to $G$ and Theorem \ref{E}
  shows that $X$ satisfies Murre's conjectures (A), (B') and (D) with
  respect to the CK decomposition $\{\Pi_k\}_{0\leq k \leq 2d}$.
  Let's now show that $X$ satisfies Murre's conjecture (B). This is
  clear because when $i \neq d$, $\Pi_i$ acts non-trivially only on
  $\CH_{\lfloor i/2 \rfloor}(X)$ because $\Pi_i$ factors through a
  curve and acts only in one degree in homology (see also \cite[Prop.
  2.9]{Vial3}); and when $i=d$, ${}^t\Pi_d = \Pi_d$ acts trivially on
  $\CH_l(X)$ for $l < n-1$ for dimension reasons (the action of
  $({}^t\Pi_d)_*$ on $\CH^l(X)$ factors through $\CH_{l+1-n}(Z_n)$ or
  the motive $(X,\Pi_d)$ is a direct summand of
  $h(\tilde{Z}_n)(n-1)$).

  The above shows that that $F^1\CH_i(X) = \CH_i(X)_\hom$ for all $i$.
  It also shows that $\CH_i(X) = (\Pi_{2i} + \Pi_{2i+1})_*\CH_i(X)$ for
  all $i \neq n-1$.  Therefore, for all $i \neq n-1$, $F^2\CH_i(X)=0$.
  (Note that homological and algebraic equivalence do not
  necessarily agree on $\CH_{n}(X)$ but do agree on $\CH_i(X)$ for
  $i\neq n-1, n$.)  For the same reasons as in the proof of Theorem
  \ref{pair}, homological and algebraic equivalence agree on
  $\CH_{n-1}(X)$. By Theorem \ref{Chowproj} (see also \cite[Prop.
  2.10]{Vial3}), we have $\ker \big(\Pi_{d-2} :
  \CH_{n-1}(X)_{\mathrm{hom}} \r \CH_{n-1}(X)_{\mathrm{hom}} \big) =
  \ker\big (AJ_{n-1} : \CH_{n-1}(X)_{\mathrm{hom}} \r J_{n-1}(X)\big)$.
  Because $\Pi_{d-1}$ acts trivially on $\CH_{n-1}(X)$, we have
  $F^2\CH_{n-1}(X) = F^3\CH_{n-1}(X)$ and then because $\Pi_d$ is the
  only remaining CK projector acting possibly non trivially on
  $\CH_{n-1}(X)$, we have $F^4\CH_{n-1}(X)=0$. We therefore see that the
  filtration does not depend on the particular choice of a special CK
  decomposition for $X$. This settles (C'') for $X$.
\end{proof}

Combining the results above gives the following.

\begin{theorem} \label{impair} Let $X$ be a smooth projective variety
  of odd dimension $d=2n+1$ with $H^{n+1}(X,\Omega_X^{n-1})=0$. If
  $\CH_0(X)_\alg, \CH_1(X)_\alg, \ldots, \CH_{n-2}(X)_\alg$ are
  representable, then $X$ satisfies Murre's conjectures (A), (B),
  (C'') and (D). \qed
\end{theorem}

\begin{examples}
 If $X$ is one of the following :
\begin{itemize}
\item a $3$-fold with $H^2(X,O_X)=0$, e.g. a Calabi--Yau $3$-fold (more
  can be said for Fano $3$-folds since they belong to $F$) or a
  complete intersection of dimension $3$;
\item a $5$-fold with $H^3(X,\Omega_X)=0$ and with representable
  $\CH_0(X)_\alg$. For instance, this could be a rationally connected
  $5$-fold with $H^3(X,\Omega_X)=0$. Examples of such varieties are
  given by blowing up a Fano complete intersection of dimension $5$
  along a smooth subvariety $Z$ satisfying $H^2(Z,O_Z)=0$. For
  instance, $Z$ could be a smooth curve contained in $X$ or, more
  interestingly, a smooth $3$-fold obtained by intersecting $X$ with
  two (high-degree) hypersurfaces.
\end{itemize}
Then $X$ satisfies Murre's conjectures (A), (B), (C'') and (D).
\end{examples}

\paragraph{4.4.3. Complete intersections of low degree.} Here is a
list of examples of complete intersections for which we can apply
Theorems \ref{pair} and \ref{impair}. All such varieties thus satisfy
Murre's conjectures (A), (B), (C'') and (D).
\begin{itemize}
\item Hypersurfaces of $\P^6$ of degree $\leq 6$ and more generally
  Fano complete intersections of dimension $5$; cf. examples above.
\item Cubic $5$- and $6$-folds : Paranjape \cite{Paranjape} and
  Koll\'ar \cite{Kollar} showed that they satisfy $\CH_0(X) = \CH_1(X)
  =\Q$. Moreover, Theorem \ref{rep} shows that the cubic $5$-fold
  satisfies Murre's conjecture (C).
\item The intersection of a quadric and of a cubic of dimension $6$ :
  Esnault, Levine and Viehweg \cite{ELV} proved that it satisfies
  $\CH_0(X) = \CH_1(X) =\Q$.
\item The intersection of two quadrics of dimension $8$ : Esnault,
  Levine and Viehweg \cite{ELV} proved that it satisfies $\CH_0(X) =
  \CH_1(X) = \CH_2(X) =\Q$.
 \item Quartic $7$-folds : Otwinowska \cite[Cor. 1]{Otwinowska} showed
  that they satisfy $\CH_0(X) = \CH_1(X) =\Q$.
\item Cubic hypersurfaces of dimensions $8$, $9$, $11$, $14$ and $15$ :
  Otwinowska \cite[Cor. 1]{Otwinowska} showed that they satisfy the
  assumptions of Theorems \ref{pair} and \ref{impair}.
\end{itemize}

More generally, Otwinowska \cite{Otwinowska} showed that if $X$ is a
smooth hyperplane section of a hypersurface in $\P^{n+1}$ covered by
$l$-planes then $\CH_i(X)_\hom = 0$ for $i \leq l-1$ (see also Esnault,
Levine and Viehweg \cite{ELV}).  Therefore when $l=\lfloor n/2\rfloor
-1$, Theorems \ref{pair} and \ref{impair} give the following theorem.

\begin{theorem} Let $l=\lfloor n/2\rfloor -1$ and let $X$ be a smooth
  hyperplane section of a hypersurface in $\P^{n+1}$ covered by
  $l$-planes. Then, $X$ satisfies Murre's conjectures (A), (B), (C'')
  and (D). \qed
\end{theorem}

\paragraph{4.4.4. Product of a surface with a variety with
  representable Chow groups.} Consider a surface $S$ and a smooth
projective variety $X$ whose Chow groups are all representable (see \S
4.4.1. for example). Thanks to \cite[Th. 4]{Vial3}, it is easy to see
that $ S \times X$ belongs to $G$. Let $Y$ be the variety obtained by
successively blowing up $X \times S$ along smooth surfaces. Then $Y$
also belongs to $G$ by Proposition \ref{belongG}.
The following is a corollary of Theorem \ref{E}.

\begin{theorem} \label{surfvar} The variety $Y$ satisfies Murre's
  conjectures (A), (B) and (D).
\end{theorem}
\begin{proof}
  By Theorem \ref{E}, there is only Murre's conjecture (B) left to
  prove. It is even enough to prove that $\Pi_i$ acts trivially on
  $\CH_l(Y)$ for $l < \lfloor i/2 \rfloor -1$. By the K\"unneth formula
  for CK decompositions and by the blow-up formula, we see that the
  idempotent $\Pi_i$ factors through a surface if $i$ is even and
  through a $3$-fold if $i$ is odd. In particular, the action of
  $\Pi_i$ on $\CH_l(Y)$ factors through $(l-\lfloor i/2 \rfloor +
  1)$-cycles on the surface or the $3$-fold.
\end{proof}

\begin{remark}
  If $Y$ is taken to be the product of a curve with a surface,
  conjectures (A), (B) and (D) for $Y$ were settled by Murre
  \cite{Murre1}.
\end{remark}

In the case when $S$ is Kimura finite-dimensional, Theorem
\ref{Beauville2} gives the following theorem.

\begin{theorem}
  If the surface $S$ is Kimura finite-dimensional and if $Z$ is
  dominated by $X \times S$, then $Z$ satisfies Murre's conjectures
  (A), (B), (C') and (D).\qed
\end{theorem}

 \begin{footnotesize}
  \bibliographystyle{plain} 
  \bibliography{bib} \medskip
 
  \noindent \textsc{DPMMS, University of Cambridge, Wilberforce Road,
    Cambridge, CB3 0WB, UK}
\end{footnotesize} 

\noindent \textit{e-mail :}  \texttt{C.Vial@dpmms.cam.ac.uk}

\end{document}